\documentclass[12pt]{article}

\usepackage{amsfonts}
\usepackage{graphicx}
\usepackage[latin1]{inputenc}
\usepackage{amsmath,amssymb,amsthm}
\usepackage{graphicx}
\usepackage{pstricks}
\usepackage{pst-tree}

\newtheorem{assumption}{Assumption}
\newtheorem{remark}{Remark}
\newtheorem{theorem}{Theorem}
\newtheorem{lemma}{Lemma}
\newtheorem{acknowledgements}{Acknowledgements}

\newcommand{\ord}{\text{ord}}
\newcommand{\nd}{\text{nd}}
\newcommand{\acn}{\text{acn}}
\newcommand{\tr}{\textbf{t}}
\newcommand{\w}{\textbf{w}}
\newcommand{\rt}{\text{root}}

\title{Taylor expansions of solutions of stochastic partial differential equations}

\author{Arnulf Jentzen\thanks{
Institute of Mathematics,
Johann Wolfgang Goethe-University,
D-60054 Frankfurt am Main, Germany,
E-mail: jentzen@math.uni-frankfurt.de}}

\begin{document}

\maketitle
\begin{abstract}   
{}
The solutions of parabolic and hyperbolic 
stochastic partial differential equations 
(SPDEs) driven by an infinite dimensional 
Brownian motion, which is a martingale, 
are in general not semi-martingales any 
more and therefore do not satisfy an 
It\^o formula like the solutions of 
finite dimensional  stochastic differential 
equations (SODEs). In particular, it is 
not possible to derive stochastic Taylor 
expansions as for the solutions of SODEs 
using an iterated application of the It\^o 
formula. However, in this article we introduce 
Taylor expansions of solutions of SPDEs 
via an alternative approach, which
avoids the need of an It\^o formula.
The main idea behind these Taylor expansions
is to use first classical Taylor expansions for 
the nonlinear coefficients of the SPDE and 
then to insert recursively the mild presentation 
of the solution of the SPDE. The iteration of 
this idea allows us to derive stochastic Taylor 
expansions of arbitrarily high order.
Combinatorial concepts of trees and woods provide 
a compact formulation of the Taylor expansions.
\end{abstract}

Taylor expansions,
stochastic partial differential equations,
SPDEs,
strong convergence,
stochastic trees


\section{Introduction}

Taylor expansions are a fundamental and repeatedly used means of approximation in mathematics, in particular in numerical analysis. Although numerical schemes for ordinary differential equations (ODEs) are often derived in an ad hoc manner, 
their convergence orders are determined by Taylor expansions
of the solution of an ODE.
An important component of these Taylor expansions are the 
iterated  total derivatives of the vector field 
corresponding higher derivatives of the solution, 
which are obtained via the chain rule, see 
Deuflhard \& Bornemann (2002).\\
An analogous situation holds for It\^{o} stochastic ordinary differential equations (SODEs), except, due to the less robust nature of stochastic calculus, the stochastic Taylor expansions here
are obtained via iterated application of the stochastic chain rule, 
the It\^{o} formula (see Kloeden \& Platen (1992)). 
Underlying this method is the fact that the solution of an SODE is an It\^{o} process or, more generally, a semi-martingale and, particular, 
of finite quadratic variation.
\\
Therefore, this approach fails, however,
if an SODE is driven by an additive stochastic 
process with infinite quadratic variation such as a fractional Brownian motion, because  the  It\^{o} formula is no longer  valid. 
A new method to derive Taylor expansions in such cases
was presented in Jentzen \& Kloeden (2007, 2008a). 
It uses the smoothness of the coefficients, 
but only minimal assumptions on the nature of the driving 
stochastic process. The resulting Taylor expansions 
there are thus robust with respect to assumptions concerning 
the driving stochastic process and, in particular, 
remain valid for other noise processes. 
The main idea in Jentzen \& Kloeden (2007) is to use 
classical Taylor
expansions for the coefficients of the SODEs (driven by an
additive H\"{o}lder continuous stochastic process) and then
to insert recursively lower order expansions into
that classical Taylor expansions.
\\
In the case of
stochastic partial differential equations (SPDEs)
the situation is different from the
finite dimensional setting.
In this article we consider SPDEs of the form 
\begin{equation} \label{eqspde}
  dU_t = \left[ A U_t + F( U_t ) \right] dt 
  + B( U_t )\,dW_t, \qquad
  U_0 = u_0
\end{equation}
on a Hilbert space $H$ and
where $ A $ is in general an unbounded linear
operator (for example $ A = \Delta $), 
$F$, $B$ are nonlinear continuous functions
and $W_t$ is a cylindrical Wiener process
(see section \ref{sec2} for a precise 
description of the equation above and 
the assumptions we use). 
The interesting thing here in the infinite dimensional
setting is, that although the SPDE \eqref{eqspde} is driven 
by the martingale Brownian motion, the solution process is not a semi-martingale any more (see Gradinaru et al. (2005) 
for a clear discussion of the problem) and therefore 
a general It\^{o} formula does not exist 
for its solutions, just special cases 
(see Gradinaru et al. (2005) and also 
Pr\'evot \& R\"ockner (2007)). Hence stochastic Taylor 
expansions for the solutions of the SPDE \eqref{eqspde} 
cannot be derived as in Kloeden \& Platen (1992)
for the solutions of finite dimensional SODEs. Consequently, until recently, only temporal approximations of low order have been derived for the solutions of such SPDEs (except for SPDEs with only finitely many
stochastic processes, see for 
example Gy\"ongy (2003)
or see also Bayer \& Teichmann (2008) 
for finitely many stochastic
processes and weak convergence).\\
The main idea of the Taylor expansions presented 
in this article
is to use classical Taylor expansions for 
the coefficients $ F $ and $ B $
in the mild integral equation version of 
the SPDE above
and then to insert the mild presentation
of the solution process into these 
classical Taylor expansions 
(see section 3 for a detailed 
presentation of these Taylor expansions).
In the case of SPDEs with additive noise
this approach was recently introduced 
in Jentzen \& Kloeden (2008c).
Moreover, numerical schemes based on these
Taylor expansions there 
have already been introduced in 
Jentzen \& Kloeden (2008b)
and in Jentzen (2008).
Indeed, it can be seen from a theoretical point
of view and in simulations
that these Taylor expansions lead to higher
order numerical schemes for SPDEs with
additive noise 
(see for example section 4.3 in Jentzen (2008)).
Here, we consider the case with general noise,
where the technical difficulty is to estimate
the reminder terms of the diffusion 
coefficient $ B $ inside the
infinite dimensional stochastic integral.
To sum up: with this approach 
we avoid the need for an It\^{o} formula but
nevertheless we can derive stochastic 
Taylor expansions 
of arbitrarily high 
order for the solution of the SPDE \eqref{eqspde}.
\\
The paper is organized as follows. 
In the next section, we describe precisely the SPDE that  we are considering  and state the assumptions that we require on its terms and coefficients and on the initial value.
Then, in the third section, we sketch the idea and 
notation for deriving simple Taylor expansions, 
which we develop in section four in some detail 
using combinatorial objects, specifically stochastic 
trees and woods, to derive Taylor expansions of an 
arbitrarily high order. We also provide an estimate 
for the reminder terms of the Taylor expansions 
there. (Proofs are postponed to the final section). 
These results are illustrated  with some representative 
examples in the fifth section. Numerical scheme 
based on these Taylor expansions will be 
discussed elsewhere.

\section{Setting and Assumptions} \label{sec2}

Fix $ T > 0 $ and 
let $ \left( \Omega, \mathcal{F}, \mathbb{P} \right) $ be 
a probability space with a 
normal filtration $ \mathcal{F}_t $, $ t \in [0,T]$, 
see e.g. Da Prato \& Zabczyk (1992) for details.
In addition, 
let $ \left( H, \left< \cdot, \cdot \right> \right) $ be a
separable Hilbert space with its norm 
denoted by $ \left| \cdot \right| $.
Moreover, let 
$ \left( D, \left| \cdot \right|_{D} \right) $ 
be a separable Banach space
with $ H \subset D $ continuously.
Then, we consider the SPDE \eqref{eqspde} 
in the mild integral equation form
\begin{equation} \label{eqmain}
  U_t = e^{ A t } u_0 + \int^{ t }_{ 0 } e^{ A (t-s) } F( U_s )\,ds
  +
  \int^{ t }_{ 0 } e^{ A (t-s) } B(U_s)\,dW_s
  \qquad
  \text{a.s.}
\end{equation}
on $H$, where $ W_t $, $ t \in [0,T]$, 
is a cylindrical $Q$-Wiener process  with $ Q = I $ 
with respect to $ \mathcal{F}_t $, $ t \in [0,T] $ on 
another separable Hilbert 
space $ \left( U, \left< \cdot, \cdot \right> \right) $ 
(space-time white noise)
and the objects $ A $, $ F $, $ B $ and $ u_0 $ are 
specified through the following assumptions.
Here $ L(U,D) $ denotes the space of 
all bounded linear operators
from $ U $ to $ D $.

\begin{assumption}\label{linearoperator} {\rm (Linear Operator $A$)} 
Let $ \mathcal{I} $ be a countable set.
Moreover, let $ ( \lambda_i )_{ i \in \mathcal{I} } $ be a family 
of positive real numbers 
with $ \inf_{ i \in \mathcal{I} } \lambda_i > 0 $
and let
$ ( e_i )_{ i \in \mathcal{I} } $ be an orthonormal basis of $H$.
Then, suppose that the linear 
operator $ A : D(A) \subset H \rightarrow H $ is given by
\[
  A v = \sum_{ i \in \mathcal{I} } 
  -\lambda_i \left< e_i, v \right> e_i
\] 
for all $ v \in D(A) $ with 
$ 
D(A) = \left\{ v \in H \big| 
\sum_{ i \in \mathcal{I} } \left| \lambda_i \right|^2 \left| 
\left< e_i, v \right>
\right|^2 < \infty \right\}
$.
\end{assumption}

\begin{assumption}\label{drift} {\rm (Drift $F$)}
The nonlinearity $ F : H \rightarrow H $ is 
infinitely often 
Fr\'{e}chet 
differentiable and its derivatives 
satisfy
$
  \sup_{ v \in H } \left| F^{(i)}(v) \right| 
  < \infty
$
for all $ i \in \mathbb{N} $.
\end{assumption}

Let $ D((-A)^r) $, $ r \in \mathbb{R} $, denote the 
interpolation spaces of powers of the operator $ -A $, see 
for example Sell \& You (2002) and let $ \left| \cdot \right|_{HS} $
denote the Hilbert-Schmidt norm for Hilbert-Schmidt operators
from $ U $ to $ H $.

\begin{assumption}\label{stochconv} {\rm (Diffusion $ B $)}
Suppose $ D \subset D((-A)^{-r}) $ continuously for some $ r \geq 0 $.
Moreover, let $ B : H \rightarrow L(U,D) $
be infinitely often 
Fr\'{e}chet differentiable and suppose that
$ e^{ A t } B^{(i)}(v)(w_1,\dots,w_i) $
and $ (-A)^\gamma e^{ A t } B(v) $
are Hilbert-Schmidt operators from $ U $ 
to $ H $ such that
\begin{eqnarray*}
  \left| 
    e^{ A t } B^{(i)}( v )(w_1, \dots, w_i) 
  \right|_{HS}
  & \leq &
    L_i \, 
   \left( 1 + |v| \right)
   \left| w_1 \right| \dots \left| w_i \right| 
   \, t^{(\delta - \frac{1}{2})} ,
  \\
  \left| e^{ A t } \left( B( v ) - B(w) \right) \right|_{HS}
  & \leq &
    L_0 \, \left| v - w \right| \, 
  t^{(\varepsilon - \frac{1}{2})} ,
  \\
  \left| (-A)^\gamma e^{ A t } B( v ) \right|_{HS}
  & \leq &
  L_0 \, \left( 1 + |v| \right) 
  \, t^{(\varepsilon - \frac{1}{2})}
\end{eqnarray*}
for all $ v,w,w_1,\dots,w_i \in H $, $ i \geq 0 $ and 
$ t \in (0,T] $ with 
constants $ L_0, L_1, \ldots > 0 $, $ \delta \in (0,\frac{1}{2}] $,
$ \gamma \in (0,1) $
and $ \varepsilon \in (0,\frac{1}{2}) $. 
\end{assumption}

\begin{assumption}\label{initial} {\rm (Initial value $ u_0 $)}\,\,
Let $ u_0 : \Omega \rightarrow D((-A)^\gamma) $ 
be a $ \mathcal{F}_0 $-measurable 
random variable
with the property that
$
  \mathbb{E} \left| (-A)^\gamma u_0 \right|^p < \infty
$
for every $ p \geq 1 $ and $ \gamma \in (0,1) $ given in Assumption \ref{drift}.
\end{assumption}
Later, it will be clear, that it would suffice
to postulate that only the first $ K $-derivatives
of $ F $ and $ B $ are bounded in
the sense above with $ K \in \mathbb{N} $
sufficiently high, but for simplicity
we use Assumption \ref{drift} and 
Assumption \ref{stochconv}.
Although similar assumptions are used in 
the literature on the approximation of 
this kind of SPDEs (see for example 
Assumption H1-H3 in Hausenblas (2003) or 
see also Lord \& Shardlow (2007),
Jentzen \& Kloeden (2008b) and
M\"uller-Gronbach \& Ritter (2007b)),
these assumptions here are in a way not
satisfying, since the nonlinear 
terms $ F $, $ B $ have to be global 
Lipschitz continuous and
Fr\'{e}chet differentiable on the space $ H $,
which is not so often fulfilled in applications.
However as in Kloeden \& Platen (1992), 
we present here a Taylor expansion 
with strong assumption and
then one can use localization techniques
as in G\"ongy (1998) or in 
Jentzen et al. (2008)
to show convergence under less restrictive
assumptions.
Furthermore, note that the stochastic 
heat equation
with additive 
or particularly multiplicative noise 
is a non trivial
example, which satisfies for example
these assumptions 
(see section \ref{secex}).
Note also that Assumption \ref{stochconv}
is in that way elegant that in contrast
to Theorem 7.4 and Theorem 
7.6 in Da Prato \& Zabczyk (1992) 
it combines space time white noise 
and trace class noise in one setting
(see Jentzen \& Kloeden (2008d)).
While in the case of space time white
noise one usually requires somethings on $ e^{ A t } $
and while in the case of trace class
noise one usually requires somethings
on $ B( \cdot ) $, here we use Assumption
\ref{stochconv}, which postulates somethings
on $ e^{ A t } B( \cdot ) $. 
In that way both cases are here 
handled in one setting.
\\
Now, we present some properties of the solution
and some notations and then start with the
Taylor expansions in the next section.
A solution of the SPDE \eqref{eqspde} 
(under the Assumptions \ref{linearoperator}-\ref{initial}) 
is
a predictable (with respect to $ \mathcal{F}_t $)
stochastic process 
$ U : \Omega \times [0,T] \rightarrow H $
with
$
  \sup_{ 0 \leq t \leq T } 
  \mathbb{E} \left| U_t \right|^2 
  < \infty 
$
and which satisfies
$$
  \mathbb{P}\left[ 
    U_t = e^{ A t } u_0
    + \int_0^t e^{ A ( t - s ) } F( U_s ) ds
    + \int_0^t e^{ A ( t - s ) } B( U_s ) dW_s
  \right]
  = 1
$$
for all $ t \in [0,T] $.
Under the Assumptions \ref{linearoperator}-\ref{initial}
there is an up to modifications unique stochastic
process $ U : \Omega \times [0,T] \rightarrow H $,
which is a solution of the SPDE \eqref{eqspde}
in the sense above and furthermore satisfies
\begin{equation} \label{boundedmom}
  \sup_{ 0 \leq t \leq T } 
  \left| (-A)^\gamma U_t \right|_{ L^p } < \infty
\end{equation}
for all $ p \geq 1 $ (see Jentzen \& Kloeden (2008d)), 
where $ \gamma \in (0,1) $ is given
in Assumption \ref{stochconv}
and where $\left| Z \right|_{ L^p }$ 
$:=$ 
$\left( \mathbb{E} \left| Z \right|^p \right)^{\frac{1}{p}} $ 
is the $ L^p $-norm of 
a random variable $ Z : \Omega \rightarrow H $.
Henceforth we fix $t_0 \in [0,T)$ and 
denote by $ \mathcal{P} $ the 
set of all equivalence classes
of predictable stochastic processes 
$$ 
  X : \Omega \times [t_0,T] \rightarrow H
  \qquad
  \text{with}
  \quad
  \sup_{ t_0 \leq t \leq T } \left| X_t \right|_{ L^p } 
  < \infty 
  \quad \forall\,  p \geq 1,
$$ 
where two processes are in one equivalence class
if they are modifications of each other.
Finally, we define the finite constants
\begin{equation} \label{constants}
  K_i := \sup_{ v \in H }
  \left| F^{(i)}(v) \right| 
  \qquad
  \text{and}
  \qquad
  R_i := 
  \sup_{ 0 \leq t \leq T } 
  \left(
  \left| (-A)^\gamma U_t \right|_{ L^i }
  +
  \left| U_t \right|_{ L^i }
  \right)
\end{equation}
for all $ i \in \mathbb{N} $.

\section{Taylor expansions} \label{sectay}

In this section we present the notation and basic idea 
behind the derivation of the Taylor expansions. 
We write
\[
  \Delta U_s := U_s - U_{ t_0 }, \qquad
  \Delta s := s - t_0
\]
for  $s$ $\in$ $[t_0,T]$ $\subset$ $[0,T]$, 
thus $ \Delta U $ denotes  the stochastic
process $ \Delta U_t$, $t \in [t_0,T]$.
Here and below $ U $ is always the unique solution process
of the SPDE \eqref{eqspde}.
Note that $ \Delta U $ is in $\mathcal{P}$. Firstly,   we introduce some integral operators and  an expression relating them, and  then we show how they can be used to derive  some simple  Taylor expansions.

\subsection{Integral Operators}
Let $ j \in \left\{ 0,1,2,{1^*}, 2^* \right\}$, where 
the indices $\left\{ 0,1,2\right\}$ will label expressions containing only a constant value or  no value of the  SPDE  solution, 
while ${1^*}$ and $ 2^*$
will label certain integrals with time dependent values of the 
solutions in the integrand. 
Specifically, we define the stochastic processes
$ I^0_j $ in $ \mathcal{P} $ by 
\[
  I^0_j(s)
  :=
 \displaystyle{  \begin{cases}
 \left( e^{ A \Delta s } - I \right) U_{t_0}
    & j=0 \\[1ex]
    \int^{s}_{t_0} e^{ A (s-r) } 
      F( U_{t_0} )\,
    dr & j = 1 \\[1ex]
    \int^{s}_{t_0} e^{ A (s-r) } 
      B( U_{ t_0 } ) \,
    dW_r & j = 2 \\[1ex]
    \int^{s}_{t_0} e^{ A (s-r) } 
      F( U_{r} )\,
    dr & j = {1^*} \\[1ex] 
    \int^{s}_{t_0} e^{ A (s-r) } 
      B( U_{r} )\,
    dW_r & j = {2^*} 
  \end{cases}}
\]
for each $ s \in [t_0,T] $. 
Given  $ i \in \mathbb{N} $ and 
$ j \in \left\{ 1,2,{1^*},2^* \right\} $, we then define the 
$i$-multilinear symmetric 
mapping 
$I^i_j$ $:$ $\mathcal{P}^i$ $:=$ $\underbrace{\mathcal{P} \times 
\dots \times \mathcal{P}}_{\text{$i$-times}}$ 
$\rightarrow$ $\mathcal{P}$ with $I^i_j[ g_1, \dots, g_i ](s)$ by
\[
 \displaystyle{  \frac{ 1 }{ i! } 
    \int^{s}_{t_0} }e^{ A (s-r) } 
      F^{(i)}( U_{t_0} )\left( g_1(r), \dots, g_i(r) \right)\, 
    dr
\] 
when $ j = 1 $ respectively
\[
 \displaystyle{  \frac{ 1 }{ i! } 
    \int^{s}_{t_0} }e^{ A (s-r) } 
      B^{(i)}( U_{t_0} )
      \left( g_1(r), \dots, g_i(r) \right)\, 
    dW_r
\] 
when $j = 2$ and
 \[  
 \int^{s}_{t_0} e^{ A (s-r) } 
      \left( 
        \int^1_0 
          F^{(i)}( U_{t_0} + \theta \Delta U_r )
          \left(
            g_1(r), \ldots,
            g_i(r)  
          \right) 
            \displaystyle{ \frac{ (1-\theta)^{(i-1)} }{ (i-1)! } }\,
        d\theta
      \right)
    dr, 
\] 
when $j = {1^*}$ respectively
 \[  
 \int^{s}_{t_0} e^{ A (s-r) } 
      \left( 
        \int^1_0 
          B^{(i)}( U_{t_0} + \theta \Delta U_r )
          \left(
            g_1(r), \ldots,
            g_i(r)  
          \right) 
            \displaystyle{ \frac{ (1-\theta)^{(i-1)} }{ (i-1)! } }\,
        d\theta
      \right)
    dW_r, 
\]
when $ j=2^*$
for all $ s \in [t_0,T] $ and $ g_1, \dots, g_i $ in $ \mathcal{P} $.
Due to Assumptions \ref{linearoperator}-\ref{initial} 
one can check that the stochastic processes
$ I^0_j $ and the mappings $ I^i_j $ are
well defined.
The proof that the processes $ I^i_j $ are predictable
(and in $ \mathcal{P} $) is a little tricky:
one shows that these processes are mean square continuous 
by using the Assumptions \ref{drift} and \ref{stochconv} and
this yields that they have a predictable version.
In the next step, we consider the mild solution of the SPDE
\eqref{eqspde}, which obviously satisfies
\begin{equation} \label{eqmaindif}
  \Delta U_s = 
  \left( e^{ A \Delta s } - I \right) U_{ t_0 } 
  + \int^{ s }_{ t_0 } e^{ A (s-r) } F( U_r ) \, dr
  +
  \int^{ s }_{ t_0 } e^{ A (s-r) } B( U_r ) \, d W_r 
\end{equation}
almost surely, or, in terms of the above  integral operators, 
\begin{equation*}
  \Delta U_s = I^0_0(s) + I^0_{{1^*}}(s) + I^0_{2^*}(s)
  \qquad
  \text{a.s.}
\end{equation*}
for $ s \in [t_0,T] $, 
which we can write symbolically in the space $ \mathcal{P}$ as
\begin{equation} \label{eqmain3}
  \Delta U = I^0_0 + I^0_{{1^*}} + I^0_{2^*}. 
\end{equation}
The stochastic processes $I^0_0 $, 
$ I^i_j[g_1, \dots, g_i] $ for $j$ $=$ $1$ or $2$
and $ i \geq 0 $
only depend on the solution at time $ t = t_0 $. These terms are therefore useful approximations
for the solution $ U_s $. 
However, the stochastic processes  
$ I^i_{{1^*}}[g_1, \dots, g_i] $ 
and
$ I^i_{{2^*}}[g_1, \dots, g_i] $
with $ i \geq 0 $
depends on the solution path $ U_s $ with $ s \in [t_0,T] $. 
In that sense the star $ * $ at the number $ 1 $ and
$ 2 $ denotes that we need a further expansion 
for these processes.
For this we will use the important formulas
\begin{eqnarray} \label{formel1}
  I^0_{{1^*}} 
  & = & I^0_1 
  + I^1_{{1^*}}[ \Delta U ] \nonumber \\[2ex]
  &=&  I^0_1
  + I^1_{{1^*}}[ I^0_0 ] 
  + I^1_{{1^*}}[ I^0_{{1^*}} ] 
  + I^1_{{1^*}}[ I^0_{2^*} ], 
\end{eqnarray}
\begin{eqnarray} \label{formel1b}
  I^0_{{2^*}} 
  & = & I^0_1
  + I^1_{{1^*}}[ \Delta U ] \nonumber \\[2ex]
  &=&  I^0_2
  + I^1_{{2^*}}[ I^0_0 ] 
  + I^1_{{2^*}}[ I^0_{{1^*}} ] 
  + I^1_{{2^*}}[ I^0_{2^*} ], 
\end{eqnarray}
which are an immediate consequence of integration by parts
and equation \eqref{eqmain3}, 
and more generally the iterated formulas
\begin{eqnarray} \label{formel2}
  I^i_{{1^*}}[ g_1, \dots, g_i ] & = & 
  I^i_1[ g_1, \dots, g_i ]
  + I^{i+1}_{{1^*}}[ \Delta U , g_1, \dots, g_i ] 
  \nonumber \\[2ex]
  & = & 
  I^i_1[ g_1, \dots, g_i ] 
  + I^{i+1}_{{1^*}}[ I^0_0 , g_1, \dots, g_i ] 
  \nonumber \\[2ex]
  &&
  + I^{i+1}_{{1^*}}[ I^0_{{1^*}} , g_1, \dots, g_i ]
  + I^{i+1}_{{1^*}}[ I^0_{2^*} , g_1, \dots, g_i ]
\end{eqnarray}
\begin{eqnarray} \label{formel2b}
  I^i_{{2^*}}[ g_1, \dots, g_i ] & = & 
  I^i_2[ g_1, \dots, g_i ]
  + I^{i+1}_{{2^*}}[ \Delta U , g_1, \dots, g_i ] 
  \nonumber \\[2ex]
  & = & 
  I^i_2[ g_1, \dots, g_i ] 
  + I^{i+1}_{{2^*}}[ I^0_0 , g_1, \dots, g_i ] 
  \nonumber \\[2ex]
  &&
  + I^{i+1}_{{2^*}}[ I^0_{{1^*}} , g_1, \dots, g_i ]
  + I^{i+1}_{{2^*}}[ I^0_{2^*} , g_1, \dots, g_i ]
\end{eqnarray}
for $g_1, \dots, g_i$ in $\mathcal{P}$, $i\geq 1$.

\subsection{Derivation of simple Taylor expansions} 
\label{secsimple}

To derive a further expansion of 
equation \eqref{eqmain3} 
we insert formula \eqref{formel1} to
the stochastic process $I^0_{{1^*}}$
and we insert formula \eqref{formel1b}
to the stochastic process $I^0_{2^*}$, 
i.e., 
\[
  I^0_{{1^*}} = I^0_1 + I^1_{{1^*}}[ I^0_0 ] 
  + I^1_{{1^*}}[ I^0_{{1^*}} ] + I^1_{{1^*}}[ I^0_{2^*} ],
\]
$$
  I^0_{{2^*}} 
  = I^0_2
  + I^1_{{2^*}}[ I^0_0 ] 
  + I^1_{{2^*}}[ I^0_{{1^*}} ] 
  + I^1_{{2^*}}[ I^0_{2^*} ]
$$
into equation \eqref{eqmain3} to obtain 
\begin{eqnarray*}
  \Delta U &=& I^0_0 + 
  \left( 
    I^0_1 + I^1_{{1^*}}[ I^0_0 ] 
    + I^1_{{1^*}}[ I^0_{{1^*}} ] + I^1_{{1^*}}[ I^0_{2^*} ]
  \right) \\
  &&
  + \left( I^0_2
  + I^1_{{2^*}}[ I^0_0 ] 
  + I^1_{{2^*}}[ I^0_{{1^*}} ] 
  + I^1_{{2^*}}[ I^0_{2^*} ] \right),
\end{eqnarray*}
which can also be written as
\begin{eqnarray} \label{tay1rest}
  \Delta U = I^0_0 + 
    I^0_1 + I^0_2 + R
\end{eqnarray}
with
$$
  R =
  I^1_{{1^*}}[ I^0_0 ] 
  + I^1_{{1^*}}[ I^0_{{1^*}} ] + I^1_{{1^*}}[ I^0_{2^*} ] .
  + I^1_{{2^*}}[ I^0_0 ] 
  + I^1_{{2^*}}[ I^0_{{1^*}} ] 
  + I^1_{{2^*}}[ I^0_{2^*} ] .
$$
If we can show that the double integral terms 
$ I^l_{j}[ I^0_k ] $ in $ R $
are sufficient small (indeed,  this will be done
in the next section), then we obtain the
approximation
\begin{equation} \label{tay1restb}
  \Delta U \approx I^0_0 + I^0_1 + I^0_2,
\end{equation}
or, using the definition of the stochastic processes $ I^0_j $,
$$
  \Delta U_t \approx 
    \left( e^{ A \Delta t } - I \right) U_{t_0}
    + \int^{t}_{t_0} e^{ A (t-r) } 
      F( U_{t_0} )\,
    dr 
    +
    \int^{t}_{t_0} e^{ A (t-r) } 
      B( U_{t_0} )\,
    dW_r
$$
for $ t \in [t_0,T] $.
Hence
\begin{equation} \label{tay1}
  U_t \approx 
    e^{ A \Delta t } U_{t_0}
    + 
    A^{-1}\left(  e^{ A \Delta t }-I\right) F( U_{t_0} )
    +
    \int^{t}_{t_0} e^{ A (t-r) } B( U_{ t_0 } )\,dW_r, 
    \; t \in [t_0,T]
\end{equation}
is an approximation for the solution of SPDE \eqref{eqspde}. Since the right hand side of equation \eqref{tay1}   depends  on the solution only at time $ t_0 $,  
it is the first non trivial Taylor expansion of the solution of the  SPDE \eqref{eqspde}. 
The reminder terms $ I^l_{j}[ I^0_k ] $ in $ R $ of this 
approximation can be estimated by
$$
  \left|
    R(t)
  \right|_{ L^2 }
  \leq C ( \Delta t )^{ ( \delta + \min( \gamma, \delta ) ) }
$$
for $ t \in [t_0,T] $ and 
with a constant $C \geq 0$ (see Theorem \ref{mainthm} in the next section).
We write
$
  Y_t = O( (\Delta t)^{ r } )
$ 
with $ r > 0 $ for a stochastic process $ Y $ in $ \mathcal{P} $ if
$
  \left| Y_t \right|_{ L^2 }
  \leq C ( \Delta t )^r
$
holds for all $ t \in [t_0,T] $ with a constant $ C > 0 $. Therefore, we have
$$
  U_t - 
    \left( 
    e^{ A \Delta t } U_{t_0}
    + A^{-1}\left(  e^{ A \Delta t }-I\right) 
    F( U_{t_0} )
    +
    \int^{t}_{t_0} e^{ A (t-r) } 
      B(U_{t_0})\,
    dW_r \right) 
    = 
    O( 
      (\Delta t)^{ ( \delta + \min( \gamma, \delta ) ) } 
    ),
$$
or 
\begin{equation} \label{onestep}
  U_t = 
    e^{ A \Delta t } U_{t_0}
    + 
  A^{-1}\left(  e^{ A \Delta t }-I\right) F( U_{t_0} )
    +
    \int^{t}_{t_0} e^{ A (t-r) } 
      B(U_{t_0})\,
    dW_r 
  + 
  O( (\Delta t)^{ ( \delta + \min( \gamma, \delta ) ) } ) .
\end{equation}
The approximation \eqref{tay1} thus has 
order $\delta + \min(\gamma,\delta)$ in 
the above strong sense.
It plays an analogous role to the simplest  strong Taylor expansion in 
Kloeden \& Platen (1992) on which 
the Euler-Maruyama scheme for finite 
dimensional SODEs is based and was in 
fact in the case of additive noise 
used in Jentzen \& Kloeden (2008b) 
to derive the exponential Euler 
scheme for the SPDE \eqref{eqspde}. 

\subsection{Higher order Taylor expansions} \label{sechigher}

Further expansions of the reminder terms 
in a Taylor expansion give a Taylor expansion 
of higher order.
To illustrate this, we will expand 
the terms 
$
  I^1_{{2^*}}[ I^0_0 ] 
$
and
$
  I^1_{{2^*}}[ I^0_{2^*} ]
$
in $ R $ in equation \eqref{tay1rest}. 
From \eqref{formel1b} 
and \eqref{formel2b} we have
\begin{eqnarray*}
  I^1_{{2^*}}[ I^0_{0} ] & = & 
  I^1_{2}[ I^0_{0} ] 
  + I^2_{{2^*}}[ I^0_0, I^0_{0} ]
  + I^2_{{2^*}}[ I^0_{{1^*}}, I^0_{0} ]
  + I^2_{{2^*}}[ I^0_{2^*}, I^0_{0} ] 
\end{eqnarray*}
and
\begin{eqnarray*}
  I^1_{{2^*}}[ I^0_{2^*} ] 
  & = & 
  I^1_{2}[ I^0_{2^*} ] 
  + I^2_{{2^*}}[ I^0_0, I^0_{2^*} ]
  + I^2_{{2^*}}[ I^0_{{1^*}}, I^0_{2^*} ]
  + I^2_{{2^*}}[ I^0_{2^*}, I^0_{2^*} ]
  \\
  & = & 
  I^1_{2}[ I^0_{2} ]
  +
  I^1_{2}[ I^1_{2^*}[ I^0_0 ] ]
  +
  I^1_{2}[ I^1_{2^*}[ I^0_{1^*} ] ]
  +
  I^1_{2}[ I^1_{2^*}[ I^0_{2^*} ] ]
  \\
  && 
  + I^2_{{2^*}}[ I^0_0, I^0_{2^*} ]
  + I^2_{{2^*}}[ I^0_{{1^*}}, I^0_{2^*} ]
  + I^2_{{2^*}}[ I^0_{2^*}, I^0_{2^*} ] ,
\end{eqnarray*}
which we insert  into equation  \eqref{tay1rest} to obtain
\begin{eqnarray*}
  \Delta U = I^0_0 + I^0_1 + I^0_2 
  + I^1_2[ I^0_0 ] + I^1_2[ I^0_2 ] + R, 
\end{eqnarray*}
where the reminder term $ R $ is here given by
\begin{multline*}
  R = 
  I^1_{{1^*}}[ I^0_0 ] 
  + I^1_{{1^*}}[ I^0_{{1^*}} ] + I^1_{{1^*}}[ I^0_{2^*} ]
  + I^1_{{2^*}}[ I^0_{{1^*}} ] 
  + I^2_{{2^*}}[ I^0_0, I^0_{0} ]
  + I^2_{{2^*}}[ I^0_{{1^*}}, I^0_{0} ]
  + I^2_{{2^*}}[ I^0_{2^*}, I^0_{0} ] 
  \\
  +
  I^1_{2}[ I^1_{2^*}[ I^0_0 ] ]
  +
  I^1_{2}[ I^1_{2^*}[ I^0_{1^*} ] ]
  +
  I^1_{2}[ I^1_{2^*}[ I^0_{2^*} ] ]
  + I^2_{{2^*}}[ I^0_0, I^0_{2^*} ]
  + I^2_{{2^*}}[ I^0_{{1^*}}, I^0_{2^*} ]
  + I^2_{{2^*}}[ I^0_{2^*}, I^0_{2^*} ] .
\end{multline*}
From Theorem \ref{mainthm} in the 
next section we will see that
$R(t)$ $=$ $O( (\Delta t)^{ 
( \delta + 2 \min( \gamma, \delta ) ) 
} ) $.  
Thus we have 
$$
  \Delta U = I^0_0 + I^0_1 + I^0_2 
  + I^1_2[ I^0_0 ] + I^1_2[ I^0_2 ] 
  + O( (\Delta t)^{ 
( \delta + 2 \min( \gamma, \delta ) ) 
} ),
$$
which can also be written as
\begin{eqnarray*}
  U_t &=& 
    e^{ A \Delta t } U_{t_0}
    + A^{-1}\left(  e^{ A \Delta t }-I\right) F( U_{t_0} )
    +
    \int^{t}_{t_0} e^{ A (t-s) } B( U_{t_0} )
    dW_s \\
  &&
  + 
  \int^t_{t_0} e^{ A ( t-s) }
  B'( U_{t_0} ) 
  \left( e^{ A \Delta s } - I \right) U_{ t_0}
  dW_s \\
  &&
  + 
  \int^t_{t_0} e^{ A ( t-s) }
  B'( U_{t_0} ) 
    \int^s_{t_0} 
      e^{ A ( s-r) } B(U_{t_0}) 
    dW_r \,
  dW_s 
  + 
  O( 
  ( \Delta t)^{ ( \delta + 2 \min( \gamma, \delta ) ) } 
  )
\end{eqnarray*}
for $ t \in [t_0,T] $.
This approximation is of 
order $ \delta + 2 \min( \gamma, \delta ) $.
The double integral $ dW_r \, dW_s $ indicates that
this Taylor expansion is in a way the infinite dimensional
analogon of the Milstein scheme for SODEs.\\
\\
By iterating this idea we can construct further Taylor expansions.
In particular, we will show in the next section how a 
Taylor expansion with arbitrary high order can be achieved.

\section{Systematic derivation of  Taylor expansions
of arbitrary high order} \label{sectrees}

The basic mechanism for deriving a Taylor expansion for the SPDE \eqref{eqspde}  was explained in the 
previous section. We will now show  how Taylor expansions of  arbitrary high order can be  derived 
and will also estimate their reminder terms. For this we will identify the terms occurring in a Taylor expansions
by combinatorial objects, i.e. trees. It is a standard tool in numerical analysis to describe higher order 
terms in a Taylor expansion via rooted trees (see, for example, 
Butcher (1987) for ODEs and
Burrage \& Burrage (2000),
R\"o$\ss$ler (2001, 2004) for SODEs). In particular, we introduce a class of trees which is appropriate for
our situation and show 
how the trees relate to the desired Taylor expansions.

\subsection{Stochastic trees and woods}
We begin with the definition of the trees that we need, 
adapting the standard notation of the trees used 
in the Taylor 
expansion of SODEs (see for example Definition 2.3.1 in 
R\"o$\ss$ler (2001) as well as 
Burrage \& Burrage (2000),
R\"o$\ss$ler (2004)).\\

Let $ N \in \mathbb{N} $ be a natural number and let 
$$
  \textbf{t}': 
  \left\{ 2, \dots, N \right\} 
  \rightarrow \left\{ 1, \dots, N \right\}, \qquad
  \textbf{t}'': 
  \left\{ 1, \dots, N \right\} 
  \rightarrow \left\{ 0, 1, 2, {1^*}, 2^* \right\}, \qquad
$$
be two mappings with the property that
$
  \textbf{t}'(j) < j
$
for all $ j \in \left\{ 2, \dots, N \right\} $. The pair of mappings 
$ \tr$ $=$ $\left( \tr', \tr'' \right) $ is a {S-tree} (stochastic tree) of length $ N$ $=$ $l(\tr)$ 
nodes.\\

Every  {S-tree} can be represented as a graph, whose
nodes are given by the set 
$ \nd( \tr ) := \left\{ 1, \dots, N \right\} $ 
and whose arcs are described by the mapping $ \tr' $
in the sense that there is an edge from $ j $ to $ \tr'(j) $
for every node $ j \in \left\{ 2, \dots, N \right\} $.
The mapping $ \tr'' $ is an additional   labelling 
of the nodes with 
$ \tr''(j) \in \left\{ 0, 1, 2, {1^*}, 2^* \right\} $ 
indicating the type of node $ j $ for every $ j \in \nd(\tr) $.
\begin{figure}
\begin{center}
\begin{picture}(100,130)
\pstree[treemode=U, dotsize=9pt]{
  \Tdot[dotsize=9pt,dotstyle=*]~*[tnpos=l]{$1$ } 
}{
  \pstree{ \Tdot[dotstyle=square*]~*[tnpos=l]{$2$ } }{
    \Tdot[dotstyle=Bo]~*[tnpos=l]{$3$ } 
    \Tdot[dotstyle=square]~*[tnpos=l]{$5$ } 
  }
  \Tdot[dotstyle=otimes]~*[tnpos=l]{$4$ } 
}
\end{picture}
\begin{picture}(150,130)
  \pstree[treemode=U, dotsize=9pt]{
    \Tdot[dotsize=9pt,dotstyle=otimes]~*[tnpos=l]{$1$ } 
  }{
    \Tdot[dotstyle=otimes]~*[tnpos=l]{$2$ }
    \Tdot[dotstyle=Bo]~*[tnpos=l]{$3$ } 
    \pstree{ \Tdot[dotstyle=*]~*[tnpos=l]{$4$ } }{
      \Tdot[dotstyle=*]~*[tnpos=l]{$6$ } 
      \Tdot[dotstyle=square]~*[tnpos=l]{$7$ } 
    }
    \Tdot[dotstyle=square*]~*[tnpos=l]{$5$ } 
  }
\end{picture}
\caption{ \label{fig1} Two examples of stochastic trees }
\end{center}
\end{figure}
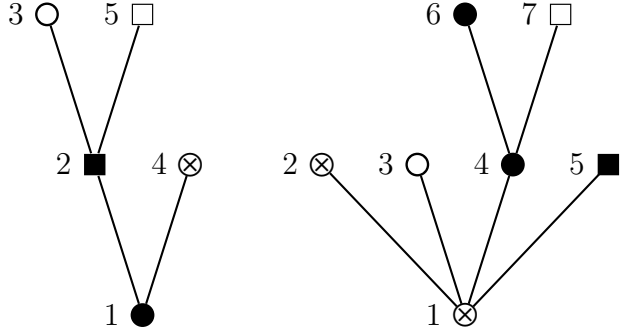
The left picture in Figure \ref{fig1} corresponds to
the tree $ \tr_1 = ( \tr_1', \tr_1'' ) $ with
$ \nd( \tr_1 ) = \left\{ 1,2,3,4,5 \right\} $ given by 
$$ 
  \tr_1'(5) = 2, \qquad
  \tr_1'(4) = 1, \qquad
  \tr_1'(3) = 2, \qquad
  \tr_1'(2)= 1
$$
and
$$
  \tr_1''(1) = 1, \qquad
  \tr_1''(2)= {1^*}, \qquad
  \tr_1''(3)= 2, \qquad 
  \tr_1''(4)= 0,
  \qquad
  \tr_1''(5)=2^* .
$$
The root is always presented as the lowest node. 
The number on the left of a node in Figure \ref{fig1} 
is the number of the node of the corresponding tree. 
The type of the nodes in Figure \ref{fig1} depends on the additional 
labelling of the nodes given by $ \tr_2'' $. 
More precisely, 
we represent a node $ j \in \nd(\tr_1) $
by 
\, \psdot*[dotsize=9pt, dotstyle=otimes](0pt,3pt) \; 
if $ \tr_1''(j) = 0 $,
by \,\psdot*[dotsize=9pt, dotstyle=*](0pt,3pt) \; 
if $ \tr_1''(j) = 1 $,
by \, \psdot*[dotsize=9pt, dotstyle=Bo](0pt,3pt) \; 
if $ \tr_1''(j) = 2 $,
by \, \psdot*[dotsize=9pt, dotstyle=square*](0pt,3pt) \; 
if $ \tr_1''(j) = 1^* $
and finally 
by \, \psdot*[dotsize=9pt, dotstyle=square](0pt,3pt) 
\; if $ \tr_1''(j) = {2^*} $. 
The right picture in Figure \ref{fig1} corresponds 
to the tree $ \tr_2 = ( \tr_2', \tr_2'' ) $ with
$ \nd( \tr_2 ) = \left\{ 1, \dots, 7 \right\} $ given
by
$$ 
  \tr_2'(7) = 4, \;
  \tr_2'(6) = 4, \;
  \tr_2'(5) = 1, \;
  \tr_2'(4) = 1, \;
  \tr_2'(3) = 1, \;
  \tr_2'(2)= 1
$$
and
$$
  \tr_2''(1) = 0, \;
  \tr_2''(2) = 0, \;
  \tr_2''(3) = 2, \;
  \tr_2''(4) = 1, \;
  \tr_2''(5) = {1^*}, \;
  \tr_2''(6)= 1, \;
  \tr_2''(7)= 2^*.
$$
We denote the set of all stochastic trees by \textbf{ST} and  will also consider a tuple of trees, i.e. a wood.
The set of  {S-woods} ({stochastic woods}) is defined by 
\[
  \textbf{SW} :=
  \bigcup^{\infty}_{l=1}
  ( \textbf{ST} )^l.
\]
Of course, we have the embedding $ \textbf{ST}$ $\subset$ $\textbf{SW} $. A simple example of an {S-wood} 
which will be required 
later is $ \w_0 $ $=$ $(\tr_1, \tr_2, \tr_3) $
with $ \tr_1 $, $ \tr_2 $ and $ \tr_3 $ are 
given by $ l(\tr_1)$ $=$ $l(\tr_2)$ $=$ $l(\tr_3)$ $=$ $1 $
and $ \tr_1''(1)$ $=$ $0 $, 
$ \tr_2''(1)$ $=$ ${1^*} $, $ \tr_3''(1)$ $=$ $2^* $. 
This is shown in
Figure \ref{fig4} where the left 
tree corresponds to $ \tr_1 $,
the middle one to $ \tr_2 $ and the right tree corresponds to $ \tr_3 $.
\begin{figure}
\begin{center}
\begin{picture}(50,15)
\pstree[treemode=U, dotsize=9pt]{ \Tdot[dotsize=9pt,dotstyle=otimes]~*[tnpos=l]{$1$ } }{}
\end{picture}
\begin{picture}(50,15)
\pstree[treemode=U, dotsize=9pt]{ \Tdot[dotsize=9pt,dotstyle=square*]~*[tnpos=l]{$1$ } }{}
\end{picture}
\begin{picture}(50,15)
\pstree[treemode=U, dotsize=9pt]{ \Tdot[dotsize=9pt,dotstyle=square]~*[tnpos=l]{$1$ } }{}
\end{picture}
\caption{ \label{fig4} The stochastic wood ${\rm \mathbf{w} }_0 $ in {\rm \textbf{SW}} }
\end{center}
\end{figure}

\subsection{Construction of stochastic trees and woods}

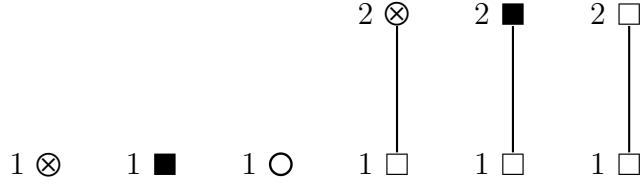
\begin{figure}
\begin{center}
\begin{picture}(40,65)
\pstree[treemode=U, dotsize=9pt]{ \Tdot[dotsize=9pt,dotstyle=otimes]~*[tnpos=l]{$1$ } }{}
\end{picture}
\begin{picture}(40,65)
\pstree[treemode=U, dotsize=9pt]{ \Tdot[dotsize=9pt,dotstyle=square*]~*[tnpos=l]{$1$ } }{}
\end{picture}
\begin{picture}(40,65)
\pstree[treemode=U, dotsize=9pt]{ \Tdot[dotsize=9pt,dotstyle=Bo]~*[tnpos=l]{$1$ } }{}
\end{picture}
\begin{picture}(40,65)
\pstree[treemode=U, dotsize=9pt]{ \Tdot[dotsize=9pt,dotstyle=square]~*[tnpos=l]{$1$ } }{
  \Tdot[dotstyle=otimes]~*[tnpos=l]{$2$ } 
}
\end{picture}
\begin{picture}(40,65)
\pstree[treemode=U, dotsize=9pt]{ \Tdot[dotsize=9pt,dotstyle=square]~*[tnpos=l]{$1$ } }{
  \Tdot[dotstyle=square*]~*[tnpos=l]{$2$ } 
}
\end{picture}
\begin{picture}(40,65)
\pstree[treemode=U, dotsize=9pt]{ \Tdot[dotsize=9pt,dotstyle=square]~*[tnpos=l]{$1$ } }{
  \Tdot[dotstyle=square]~*[tnpos=l]{$2$ } 
}
\end{picture}
\caption{ \label{fig2} The stochastic wood ${\rm \mathbf{w} }_1 $ in {\rm \textbf{SW}} }
\end{center}
\end{figure}

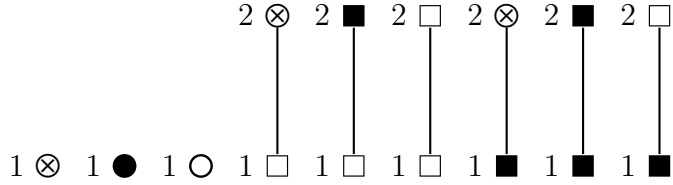
\begin{figure}
\begin{center}
\begin{picture}(25,65)
\pstree[treemode=U, dotsize=9pt]{ \Tdot[dotsize=9pt,dotstyle=otimes]~*[tnpos=l]{$1$ } }{}
\end{picture}
\begin{picture}(25,65)
\pstree[treemode=U, dotsize=9pt]{ \Tdot[dotsize=9pt,dotstyle=*]~*[tnpos=l]{$1$ } }{}
\end{picture}
\begin{picture}(25,65)
\pstree[treemode=U, dotsize=9pt]{ \Tdot[dotsize=9pt,dotstyle=Bo]~*[tnpos=l]{$1$ } }{}
\end{picture}
\begin{picture}(25,65)
\pstree[treemode=U, dotsize=9pt]{ \Tdot[dotsize=9pt,dotstyle=square]~*[tnpos=l]{$1$ } }{
  \Tdot[dotstyle=otimes]~*[tnpos=l]{$2$ } 
}
\end{picture}
\begin{picture}(25,70)
\pstree[treemode=U, dotsize=9pt]{ \Tdot[dotsize=9pt,dotstyle=square]~*[tnpos=l]{$1$ } }{
  \Tdot[dotstyle=square*]~*[tnpos=l]{$2$ } 
}
\end{picture}
\begin{picture}(25,70)
\pstree[treemode=U, dotsize=9pt]{ \Tdot[dotsize=9pt,dotstyle=square]~*[tnpos=l]{$1$ } }{
  \Tdot[dotstyle=square]~*[tnpos=l]{$2$ } 
}
\end{picture}
\begin{picture}(25,65)
\pstree[treemode=U, dotsize=9pt]{ \Tdot[dotsize=9pt,dotstyle=square*]~*[tnpos=l]{$1$ } }{
  \Tdot[dotstyle=otimes]~*[tnpos=l]{$2$ } 
}
\end{picture}
\begin{picture}(25,70)
\pstree[treemode=U, dotsize=9pt]{ \Tdot[dotsize=9pt,dotstyle=square*]~*[tnpos=l]{$1$ } }{
  \Tdot[dotstyle=square*]~*[tnpos=l]{$2$ } 
}
\end{picture}
\begin{picture}(25,70)
\pstree[treemode=U, dotsize=9pt]{ \Tdot[dotsize=9pt,dotstyle=square*]~*[tnpos=l]{$1$ } }{
  \Tdot[dotstyle=square]~*[tnpos=l]{$2$ } 
}
\end{picture}
\caption{ \label{fig5} The stochastic wood ${\rm \mathbf{w}}_2 $ in {\rm \textbf{SW}} }
\end{center}
\end{figure}

We  define an operator  on the set \textbf{SW}, that will enable us to construct an appropriate  stochastic wood step by step.  Let $ \w$ $=$ $( \tr_1, \dots, \tr_l ) $ be a {S-wood} with $ \tr_i$ $=$ $( \tr_i', \tr_i'' ) $  for $1$ $\leq$ $i$ $\leq$ $l$. 
Moreover, let $ i \in \left\{ 1, \dots, l \right\} $ and $ j \in \left\{ 1, \dots, l(\tr_i) \right\} $ be given and 
suppose  
that $ \tr_i''( j )$ $=$ ${1^*}$ or 
$ \tr_i''( j ) $ $=$ $2^* $, 
in which case  we call the pair $ (i,j) $ an active node of $ \w $. We denote  the set of all active nodes of $ \w $ by $\acn( \w )$.
In Figures for woods respectively trees
(for example Figure \ref{fig1})
these nodes are represented by a square 
(a filled 
square 
\, \psdot*[dotsize=9pt, dotstyle=square*](0pt,3pt) 
\; \ for $1^*$ and a simple 
square \, \psdot*[dotsize=9pt, dotstyle=square](0pt,3pt) \; \ for 
$2^*$).
Now, we introduce the trees $ \tr_{l+1}$ $=$ $( \tr_{l+1}', \tr_{l+1}'' ) $, $ \tr_{l+2}$ $=$ $( \tr_{l+2}', \tr_{l+2}'' ) $ and $ \tr_{l+3}$ $=$ $( \tr_{l+3}', \tr_{l+3}'' ) $ in \textbf{ST} by
$$
\nd( \tr_{l+k} )  = \left\{ 1, \dots, l(\tr_i), l(\tr_i)+1 \right\} ,
$$
$$
\tr_{l+k}'( n )  =  \tr_i'( n ), \;\,  2 \leq n \leq l(\tr_i), \qquad  \tr_{l+k}''( n )  = \tr_i''( n ), \;\,  1
 \leq n \leq l( \tr_i ) ,
$$
$$
  \tr_{l+k}'\left( l(\tr_i)+1 \right) = j, \qquad
  \tr_{l+k}''\left( l(\tr_i)+1 \right) = 
  \begin{cases}
    0 & k = 1 \\
    1^* & k=2 \\
    2^* & k=3 \; ,
  \end{cases} 
$$
for $ k = 1,2,3 $.
Finally, we consider 
the {S-tree} 
$ \tilde{\tr}$ $=$ $\left( \tilde{\tr}', \tilde{\tr}'' \right) $
given by  $ \tilde{ \tr }'$ $=$ $\tr_i' $, but 
with $ \tilde{ \tr }'' $ given by
$
  \tilde{ \tr }''(n) = \tr_i''(n) 
$ for $ n \neq j $
and
$
  \tilde{ \tr }''(j) =
  \begin{cases}
    1 & \tr_i''(j) = 1^* \\
    2 & \tr_i''(j) = 2^*
  \end{cases} .
$
Then, we define 
$$ 
  E_{(i,j)}( \tr_1, \dots, \tr_l )
  :=
  ( \tr_1, \dots, \tr_{i-1}, \tilde{\tr}, 
  \tr_{i+1}, \dots \tr_{l+3} ) 
$$
and  consider the set of all woods that can be constructed by iteratively 
applying the $ E_{(i,j)} $ operations, i.e.
we define
$$
  \textbf{SW'} 
  :=
  \left\{
    \w \in \textbf{SW} \; \bigg|
    \begin{array}{l}
      \w = E_{(i_n,j_n)} \dots E_{(i_1,j_1)} \w_0, \,\, 
      \, n \geq 0 \\
      ( i_l, j_l ) \in \acn( E_{ (i_{l-1}, j_{l-1}) } \ldots 
      E_{(i_1,j_1)} \w_{0} ), 1 \leq l \leq n
    \end{array}
  \right\} 
$$
for the $\w_0$ introduced above. To illustrate these 
definitions, we present some examples 
using the initial stochastic wood $ \w_0 $
given in Figure \ref{fig4}. First, 
the active nodes of $ \w_0 $ 
are $\acn( \w_0 )$ $=$ $\left\{ (2,1), (3,1) \right\}$, 
since the first node in the second tree 
and the first node in the third tree are 
represented by squares.
Hence,  $ E_{(3,1)} \w_0$ 
is well defined and the resulting stochastic 
wood $ \w_1$ $=$ $E_{(3,1)} \w_0 $, which has six 
trees, is presented in Figure \ref{fig2}. 
Writing $ \w_1$ $=$ $( \tr_1, \dots, \tr_6 )$, 
the left tree in Figure \ref{fig2} corresponds 
to $ \tr_1 $, the second tree 
in Figure \ref{fig2} corresponds to $ \tr_2 $, and so on.
Moreover, we have
\begin{equation} \label{anodefig2}
  \acn( \w_1 ) =
  \left\{ 
    (2,1), (4,1), (5,1), (5,2), (6,1), (6,2)
  \right\}
\end{equation}
for the active nodes of $ \w_1 $, so 
$ \w_2$ $=$ $E_{(2,1)} \w_1 $ is 
also well defined. 
It is presented in Figure \ref{fig5}.
In Figure \ref{fig6}, we present the
stochastic wood $ \w_3 = E_{ (4,1) } \w_2 $, which
is well defined since 
\begin{equation} \label{anodefig5}
  \acn( \w_2 ) =
  \left\{
    (4,1),(5,1), (5,2), (6,1), (6,2), 
    (7,1), (8,1), (8,2), (9,1),	 (9,2) 
  \right\} .
\end{equation}
\begin{figure}
\begin{center}
\begin{picture}(25,65)
\pstree[treemode=U, dotsize=9pt]{ \Tdot[dotsize=9pt,dotstyle=otimes]~*[tnpos=l]{$1$ } }{}
\end{picture}
\begin{picture}(25,65)
\pstree[treemode=U, dotsize=9pt]{ \Tdot[dotsize=9pt,dotstyle=*]~*[tnpos=l]{$1$ } }{}
\end{picture}
\begin{picture}(25,65)
\pstree[treemode=U, dotsize=9pt]{ \Tdot[dotsize=9pt,dotstyle=Bo]~*[tnpos=l]{$1$ } }{}
\end{picture}
\begin{picture}(25,65)
\pstree[treemode=U, dotsize=9pt]{  
  \Tdot[dotsize=9pt,dotstyle=Bo]~*[tnpos=l]{$1$ } 
}{
  \Tdot[dotstyle=otimes]~*[tnpos=l]{$2$ } 
}
\end{picture}
\begin{picture}(25,70)
\pstree[treemode=U, dotsize=9pt]{ \Tdot[dotsize=9pt,dotstyle=square]~*[tnpos=l]{$1$ } }{
  \Tdot[dotstyle=square*]~*[tnpos=l]{$2$ } 
}
\end{picture}
\begin{picture}(25,70)
\pstree[treemode=U, dotsize=9pt]{ \Tdot[dotsize=9pt,dotstyle=square]~*[tnpos=l]{$1$ } }{
  \Tdot[dotstyle=square]~*[tnpos=l]{$2$ } 
}
\end{picture}
\begin{picture}(25,65)
\pstree[treemode=U, dotsize=9pt]{ \Tdot[dotsize=9pt,dotstyle=square*]~*[tnpos=l]{$1$ } }{
  \Tdot[dotstyle=otimes]~*[tnpos=l]{$2$ } 
}
\end{picture}
\begin{picture}(25,70)
\pstree[treemode=U, dotsize=9pt]{ \Tdot[dotsize=9pt,dotstyle=square*]~*[tnpos=l]{$1$ } }{
  \Tdot[dotstyle=square*]~*[tnpos=l]{$2$ } 
}
\end{picture}
\begin{picture}(25,70)
\pstree[treemode=U, dotsize=9pt]{ \Tdot[dotsize=9pt,dotstyle=square*]~*[tnpos=l]{$1$ } }{
  \Tdot[dotstyle=square]~*[tnpos=l]{$2$ } 
}
\end{picture}
\begin{picture}(55,70)
\pstree[treemode=U, dotsize=9pt]{ 
  \Tdot[dotsize=9pt,dotstyle=square]~*[tnpos=l]{$1$ } 
}{
  \Tdot[dotstyle=otimes]~*[tnpos=l]{$2$ } 
  \Tdot[dotstyle=otimes]~*[tnpos=l]{$3$ } 
}
\end{picture}
\begin{picture}(55,70)
\pstree[treemode=U, dotsize=9pt]{ 
  \Tdot[dotsize=9pt,dotstyle=square]~*[tnpos=l]{$1$ } 
}{
  \Tdot[dotstyle=otimes]~*[tnpos=l]{$2$ } 
  \Tdot[dotstyle=square*]~*[tnpos=l]{$3$ } 
}
\end{picture}
\begin{picture}(55,70)
\pstree[treemode=U, dotsize=9pt]{ 
  \Tdot[dotsize=9pt,dotstyle=square]~*[tnpos=l]{$1$ } 
}{
  \Tdot[dotstyle=otimes]~*[tnpos=l]{$2$ } 
  \Tdot[dotstyle=square]~*[tnpos=l]{$3$ } 
}
\end{picture}
\caption{ \label{fig6} The stochastic wood ${\rm \mathbf{w}}_3 $ 
in {\rm \textbf{SW}} }
\end{center}
\end{figure}
For the S-wood $ \w_3 $, we have
\begin{equation} \label{anodefig6}
  \acn( \w_3 ) =
  \left\{
    \begin{array}{c}
      (5,1), (5,2), (6,1), (6,2), (7,1), 
      (8,1), (8,2), \\ (9,1), (9,2), 
     (10,1), (11,1), (11,3), (12,1), (12,3)
    \end{array}
  \right\} .
\end{equation}
Therefore, we present the stochastic 
wood $ \w_4 = E_{ (6,1) } \w_3 $ in Figure \ref{fig7}.
\begin{figure}
\begin{center}
\begin{picture}(25,65)
\pstree[treemode=U, dotsize=9pt]{ \Tdot[dotsize=9pt,dotstyle=otimes]~*[tnpos=l]{$1$ } }{}
\end{picture}
\begin{picture}(25,65)
\pstree[treemode=U, dotsize=9pt]{ \Tdot[dotsize=9pt,dotstyle=*]~*[tnpos=l]{$1$ } }{}
\end{picture}
\begin{picture}(25,65)
\pstree[treemode=U, dotsize=9pt]{ \Tdot[dotsize=9pt,dotstyle=Bo]~*[tnpos=l]{$1$ } }{}
\end{picture}
\begin{picture}(25,65)
\pstree[treemode=U, dotsize=9pt]{  
  \Tdot[dotsize=9pt,dotstyle=Bo]~*[tnpos=l]{$1$ } 
}{
  \Tdot[dotstyle=otimes]~*[tnpos=l]{$2$ } 
}
\end{picture}
\begin{picture}(25,70)
\pstree[treemode=U, dotsize=9pt]{ \Tdot[dotsize=9pt,dotstyle=square]~*[tnpos=l]{$1$ } }{
  \Tdot[dotstyle=square*]~*[tnpos=l]{$2$ } 
}
\end{picture}
\begin{picture}(25,70)
\pstree[treemode=U, dotsize=9pt]{
  \Tdot[dotsize=9pt,dotstyle=Bo]~*[tnpos=l]{$1$ } 
}{
  \Tdot[dotstyle=square]~*[tnpos=l]{$2$ } 
}
\end{picture}
\begin{picture}(25,65)
\pstree[treemode=U, dotsize=9pt]{ \Tdot[dotsize=9pt,dotstyle=square*]~*[tnpos=l]{$1$ } }{
  \Tdot[dotstyle=otimes]~*[tnpos=l]{$2$ } 
}
\end{picture}
\begin{picture}(25,70)
\pstree[treemode=U, dotsize=9pt]{ \Tdot[dotsize=9pt,dotstyle=square*]~*[tnpos=l]{$1$ } }{
  \Tdot[dotstyle=square*]~*[tnpos=l]{$2$ } 
}
\end{picture}
\begin{picture}(25,70)
\pstree[treemode=U, dotsize=9pt]{ \Tdot[dotsize=9pt,dotstyle=square*]~*[tnpos=l]{$1$ } }{
  \Tdot[dotstyle=square]~*[tnpos=l]{$2$ } 
}
\end{picture}
\begin{picture}(55,70)
\pstree[treemode=U, dotsize=9pt]{ 
  \Tdot[dotsize=9pt,dotstyle=square]~*[tnpos=l]{$1$ } 
}{
  \Tdot[dotstyle=otimes]~*[tnpos=l]{$2$ } 
  \Tdot[dotstyle=otimes]~*[tnpos=l]{$3$ } 
}
\end{picture}
\begin{picture}(55,70)
\pstree[treemode=U, dotsize=9pt]{ 
  \Tdot[dotsize=9pt,dotstyle=square]~*[tnpos=l]{$1$ } 
}{
  \Tdot[dotstyle=otimes]~*[tnpos=l]{$2$ } 
  \Tdot[dotstyle=square*]~*[tnpos=l]{$3$ } 
}
\end{picture}
\begin{picture}(55,70)
\pstree[treemode=U, dotsize=9pt]{ 
  \Tdot[dotsize=9pt,dotstyle=square]~*[tnpos=l]{$1$ } 
}{
  \Tdot[dotstyle=otimes]~*[tnpos=l]{$2$ } 
  \Tdot[dotstyle=square]~*[tnpos=l]{$3$ } 
}
\end{picture}
\begin{picture}(55,70)
\pstree[treemode=U, dotsize=9pt]{
  \Tdot[dotsize=9pt,dotstyle=square]~*[tnpos=l]{$1$ } 
}{
  \Tdot[dotstyle=square]~*[tnpos=l]{$2$ } 
  \Tdot[dotstyle=otimes]~*[tnpos=l]{$3$ } 
}
\end{picture}
\begin{picture}(55,70)
\pstree[treemode=U, dotsize=9pt]{
  \Tdot[dotsize=9pt,dotstyle=square]~*[tnpos=l]{$1$ } 
}{
  \Tdot[dotstyle=square]~*[tnpos=l]{$2$ } 
  \Tdot[dotstyle=square*]~*[tnpos=l]{$3$ } 
}
\end{picture}
\begin{picture}(55,70)
\pstree[treemode=U, dotsize=9pt]{
  \Tdot[dotsize=9pt,dotstyle=square]~*[tnpos=l]{$1$ } 
}{
  \Tdot[dotstyle=square]~*[tnpos=l]{$2$ } 
  \Tdot[dotstyle=square]~*[tnpos=l]{$3$ } 
}
\end{picture}
\caption{ \label{fig7} The stochastic 
wood ${\rm \mathbf{w}}_4 $ in {\rm \textbf{SW}} }
\end{center}
\end{figure}
Hence, we obtain
\begin{equation} \label{anodefig7}
  \acn( \w_4 ) =
  \left\{
    \begin{array}{c}
      (5,1), (5,2), (6,2),
      (7,1), (8,1), (8,2), (9,1), (9,2), \\ (10,1),
      (11,1), (11,3), (12,1), (12,3), (13,1), \\ (13,2),
      (14,1), (14,2), (14,3), (15,1), (15,2), (15,3) 
    \end{array}
  \right\} .
\end{equation}
Finally, we present the stochastic wood
$ \w_5 = E_{(6,2)} \w_4 $
with
\begin{equation} \label{anodefig8}
  \acn( \w_5 ) =
  \left\{
    \begin{array}{c}
      (5,1), (5,2), (7,1),
      (8,1), (8,2), (9,1), (9,2), (10,1), 
      (11,1), \\ (11,3), 
      (12,1), (12,3), (13,1), (13,2),  
      (14,1), (14,2), (14,3), \\ (15,1), 
      (15,2), (15,3), 
      (16,2), (17,2), (17,3), (18,2), (18,3) 
    \end{array}
  \right\}
\end{equation}
in Figure \ref{fig8}.
\begin{figure}
\begin{center}
\begin{picture}(25,65)
\pstree[treemode=U, dotsize=9pt]{ \Tdot[dotsize=9pt,dotstyle=otimes]~*[tnpos=l]{$1$ } }{}
\end{picture}
\begin{picture}(25,65)
\pstree[treemode=U, dotsize=9pt]{ \Tdot[dotsize=9pt,dotstyle=*]~*[tnpos=l]{$1$ } }{}
\end{picture}
\begin{picture}(25,65)
\pstree[treemode=U, dotsize=9pt]{ \Tdot[dotsize=9pt,dotstyle=Bo]~*[tnpos=l]{$1$ } }{}
\end{picture}
\begin{picture}(25,65)
\pstree[treemode=U, dotsize=9pt]{  
  \Tdot[dotsize=9pt,dotstyle=Bo]~*[tnpos=l]{$1$ } 
}{
  \Tdot[dotstyle=otimes]~*[tnpos=l]{$2$ } 
}
\end{picture}
\begin{picture}(25,70)
\pstree[treemode=U, dotsize=9pt]{ \Tdot[dotsize=9pt,dotstyle=square]~*[tnpos=l]{$1$ } }{
  \Tdot[dotstyle=square*]~*[tnpos=l]{$2$ } 
}
\end{picture}
\begin{picture}(25,70)
\pstree[treemode=U, dotsize=9pt]{
  \Tdot[dotsize=9pt,dotstyle=Bo]~*[tnpos=l]{$1$ } 
}{
  \Tdot[dotstyle=Bo]~*[tnpos=l]{$2$ } 
}
\end{picture}
\begin{picture}(25,65)
\pstree[treemode=U, dotsize=9pt]{ \Tdot[dotsize=9pt,dotstyle=square*]~*[tnpos=l]{$1$ } }{
  \Tdot[dotstyle=otimes]~*[tnpos=l]{$2$ } 
}
\end{picture}
\begin{picture}(25,70)
\pstree[treemode=U, dotsize=9pt]{ \Tdot[dotsize=9pt,dotstyle=square*]~*[tnpos=l]{$1$ } }{
  \Tdot[dotstyle=square*]~*[tnpos=l]{$2$ } 
}
\end{picture}
\begin{picture}(25,70)
\pstree[treemode=U, dotsize=9pt]{ \Tdot[dotsize=9pt,dotstyle=square*]~*[tnpos=l]{$1$ } }{
  \Tdot[dotstyle=square]~*[tnpos=l]{$2$ } 
}
\end{picture}
\begin{picture}(55,70)
\pstree[treemode=U, dotsize=9pt]{ 
  \Tdot[dotsize=9pt,dotstyle=square]~*[tnpos=l]{$1$ } 
}{
  \Tdot[dotstyle=otimes]~*[tnpos=l]{$2$ } 
  \Tdot[dotstyle=otimes]~*[tnpos=l]{$3$ } 
}
\end{picture}
\begin{picture}(55,70)
\pstree[treemode=U, dotsize=9pt]{ 
  \Tdot[dotsize=9pt,dotstyle=square]~*[tnpos=l]{$1$ } 
}{
  \Tdot[dotstyle=otimes]~*[tnpos=l]{$2$ } 
  \Tdot[dotstyle=square*]~*[tnpos=l]{$3$ } 
}
\end{picture}
\begin{picture}(55,70)
\pstree[treemode=U, dotsize=9pt]{ 
  \Tdot[dotsize=9pt,dotstyle=square]~*[tnpos=l]{$1$ } 
}{
  \Tdot[dotstyle=otimes]~*[tnpos=l]{$2$ } 
  \Tdot[dotstyle=square]~*[tnpos=l]{$3$ } 
}
\end{picture}
\begin{picture}(55,70)
\pstree[treemode=U, dotsize=9pt]{
  \Tdot[dotsize=9pt,dotstyle=square]~*[tnpos=l]{$1$ } 
}{
  \Tdot[dotstyle=square]~*[tnpos=l]{$2$ } 
  \Tdot[dotstyle=otimes]~*[tnpos=l]{$3$ } 
}
\end{picture}
\begin{picture}(55,70)
\pstree[treemode=U, dotsize=9pt]{
  \Tdot[dotsize=9pt,dotstyle=square]~*[tnpos=l]{$1$ } 
}{
  \Tdot[dotstyle=square]~*[tnpos=l]{$2$ } 
  \Tdot[dotstyle=square*]~*[tnpos=l]{$3$ } 
}
\end{picture}
\begin{picture}(55,70)
\pstree[treemode=U, dotsize=9pt]{
  \Tdot[dotsize=9pt,dotstyle=square]~*[tnpos=l]{$1$ } 
}{
  \Tdot[dotstyle=square]~*[tnpos=l]{$2$ } 
  \Tdot[dotstyle=square]~*[tnpos=l]{$3$ } 
}
\end{picture}
\begin{picture}(20,135)
\pstree[treemode=U, dotsize=9pt]{
  \Tdot[dotsize=9pt,dotstyle=Bo]~*[tnpos=l]{$1$ } 
}{
  \pstree[treemode=U, dotsize=9pt]{
    \Tdot[dotstyle=square]~*[tnpos=l]{$2$ } 
  }{
    \Tdot[dotstyle=otimes]~*[tnpos=l]{$3$ }
  }
}
\end{picture}
\begin{picture}(20,135)
\pstree[treemode=U, dotsize=9pt]{
  \Tdot[dotsize=9pt,dotstyle=Bo]~*[tnpos=l]{$1$ } 
}{
  \pstree[treemode=U, dotsize=9pt]{
    \Tdot[dotstyle=square]~*[tnpos=l]{$2$ } 
  }{
    \Tdot[dotstyle=square*]~*[tnpos=l]{$3$ }
  }
}
\end{picture}
\begin{picture}(20,135)
\pstree[treemode=U, dotsize=9pt]{
  \Tdot[dotsize=9pt,dotstyle=Bo]~*[tnpos=l]{$1$ } 
}{
  \pstree[treemode=U, dotsize=9pt]{
    \Tdot[dotstyle=square]~*[tnpos=l]{$2$ } 
  }{
    \Tdot[dotstyle=square]~*[tnpos=l]{$3$ }
  }
}
\end{picture}
\caption{ \label{fig8} 	The stochastic wood ${\rm \mathbf{w}}_5 $ 
in {\rm \textbf{SW}} }
\end{center}
\end{figure}
By definition the S-woods 
$ \w_0, \w_1, \dots, \w_5 $ are in \textbf{SW}', but 
the stochastic wood given in Figure \ref{fig1} is not in \textbf{SW}'.

\subsection{Subtrees}

\begin{figure}
\begin{center}
\begin{picture}(50,65)
\pstree[treemode=U, dotsize=9pt]{
  \Tdot[dotstyle=otimes,dotsize=9pt]~*[tnpos=l]{$1$ } 
}{}
\end{picture}
\begin{picture}(50,65)
\pstree[treemode=U, dotsize=9pt]{ 
  \Tdot[dotstyle=Bo,dotsize=9pt]~*[tnpos=l]{$1$ }  
}{}
\end{picture}
\begin{picture}(65,65)
\pstree[treemode=U, dotsize=9pt]{ 
  \Tdot[dotstyle=*,dotsize=9pt]~*[tnpos=l]{$1$ } 
}{
  \Tdot[dotstyle=*]~*[tnpos=l]{$2$ } 
  \Tdot[dotstyle=square]~*[tnpos=l]{$3$ } 
}
\end{picture}
\begin{picture}(50,65)
\pstree[treemode=U, dotsize=9pt]{
  \Tdot[dotstyle=square*,dotsize=9pt]~*[tnpos=l]{$1$ } 
}{}
\end{picture}
\caption{ \label{fig3} Subtrees of the 
right tree in Figure \ref{fig1} }
\end{center}
\end{figure}
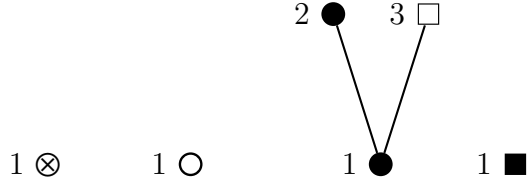

Let $ \tr = (\tr',\tr'') $ be a given S-tree with  $ \nd( \tr ) \geq 2 $.
For two nodes $ k, l \in \nd(\tr) $ with $ k \leq l $
we say that $ l $ is a grandchild of $ k $ if there
exists a sequence $ k_1$ $=$ $k$ $<$ $k_2$ $<$ $\ldots$ $<$ $k_n$ $=$ $l$ of nodes with 
$ n \geq 1$ such that
$$
  \tr'( k_{v+1} ) = k_v
  \qquad
  \text{for}
  \qquad
  v=1, \dots, n-1 .
$$
Suppose now that $ j_1 < \dots < j_n $ with $ n \geq 1 $
are the nodes of $ \tr $ such that  
$
  \tr'( j_i ) = \rt( \tr )
$
for $i$ $=$ $1$, $\dots$, $n$. 
For given
$i$ $\in$ $\left\{ 1, \dots, n \right\}$ suppose that
$$ 
  j_i = j_{i,1} < j_{i,2} < \ldots < j_{i,l_i} \leq l(\tr) ,
$$
where $ l_i \geq 1$ 
and $ j_{i,1}, \dots, j_{i,l_i} $ are 
the grandchildren of $j_i$. 
Then, we define the trees 
$ \tr_i = \left( \tr_i', \tr_i'' \right) $
with 
$
  l( \tr_i ) := l_i 
$
such that
$$
  j_{i,\tr_i'( k )} = \tr'( j_{i,k} ), \qquad
  \tr_i''( k ) = \tr''( j_{i,k} )
$$
for all $ k \in \left\{ 2, \dots, l_i \right\} $
and $ \tr_i''( 1 ) = \tr''( j_i ) $ for
every $ i = 1, \dots, n $.
We call the trees $ \tr_1, \dots, \tr_n $ defined in this way 
the \underline{subtrees} of $ \tr $.
For example, the subtrees of the right tree 
in Figure \ref{fig1}
are presented in Figure \ref{fig3}.

\subsection{Order of a tree} \label{secorder}

We also require the order of a stochastic tree 
and of a stochastic wood.
For this we introduce 
a function $\ord$ $:$ $\textbf{ST}$ $\rightarrow$ $[0,\infty)$ 
given by
\begin{multline*}
  \ord( \tr ) :=  
  \left| 
    \left\{  
      j \in \nd( \tr ) 
      \big| 
      \tr''(j) = 1 \text{ or } \tr''(j) = 1^*
    \right\} 
  \right| 
  \\
  + 
  \gamma 
  \left| 
    \left\{ 
      j \in \nd( \tr ) 
      \big| 
      \tr''(j) = 0 
    \right\} 
  \right| 
  + 
  \delta 
  \left| 
    \left\{  
      j \in \nd( \tr ) 
      \big| 
      \tr''(j) = 2 \text{ or } \tr''(j) = 2^*
    \right\} 
  \right| 
\end{multline*}
for every S-tree $ \tr = ( \tr', \tr'') $
with $ \tr''(1) \neq 0 $ and 
$ \ord(\tr) := \gamma $ 
for every S-tree $ \tr = ( \tr', \tr'') $
with $ \tr''(1) = 0 $.
For example, the order of the left 
tree in Figure \ref{fig1} is 
$ 2 + \gamma + 2 \delta $ (since
the left tree has one node 
of type $0$, two nodes of type
$1$ respectively $1^*$
and finally also two nodes
of type $2$ respectively $2^*$)
and the order of the right tree in Figure \ref{fig1} is 
$ \gamma $ since its root is of type $0$.
In addition, we say that a tree $ \tr = ( \tr', \tr'') $ 
in \textbf{ST} is \underline{active}
if there is a $ j \in \nd( \tr ) $ such 
that $ \tr''(j) = {1^*} $
or $ \tr''(j) = 2^* $.
In that sense a S-tree is active if it has an active node. 
Moreover, we define 
the order of an 
S-wood $ \w$ $=$ $( \tr_1, \dots, \tr_n )$ with $ n \geq 1 $
as 
$$
  \ord( \w ) := 
  \min
  \left\{ \ord( \tr_i ), 1 \leq i \leq n \big|
  \tr_i \text{ is active}
  \right\} .
$$
To illustrate this definition, 
we calculate the order of some 
stochastic woods.
First of all, the stochastic wood 
in Figure \ref{fig4} has order
$ \delta $, since the middle tree 
and the last (the third) tree 
in Figure \ref{fig4} are active.
More precisely, the nodes $(2,1)$ and $(3,1)$ 
of the S-wood $ \w_0 $ are active nodes
and therefore the second and the 
third tree are active.
The second tree in Figure \ref{fig4} has order
$ 1 $ (since it only consists of one node of type $ 1^* $)
and the third three in Figure \ref{fig4} has order
$ \delta $.
Hence, since $ \delta \leq 1 $, 
the S-wood $ \w_0 $ has order $ \delta $.
Since the second tree and the last three trees 
are active in the stochastic 
wood $ \w_1 $ in Figure \ref{fig2} (see equation \eqref{anodefig2} 
for the active nodes of $\w_1 $), we obtain
that the stochastic wood in Figure \ref{fig2} has order 
$ \delta + \min(\gamma,\delta) $.
The second tree in $ \w_1 $ has order $ 1 $ and
the last three trees in the S-wood $ \w_1 $ have order
$ \delta + \gamma $, $ \delta + 1 $ respectively 
$ 2\delta $.
As a third example, we consider the S-wood $ \w_2 $ 
in Figure \ref{fig5}.
The active nodes of $ \w_2 $ 
are presented in equation \eqref{anodefig5}.
Hence, the last six S-trees are active.
They have the orders $ \delta + \gamma $, 
$ \delta + 1 $, $ 2 \delta $, 
$ 1 + \gamma $, $ 2 $ and $ 1 + \delta $.
The minimum of the six real numbers is 
$ \delta + \min(\delta,\gamma) $.
Therefore, the order of the S-wood $ \w_2 $ in Figure \ref{fig5}
is also $ \delta + \min(\delta, \gamma) $.
A similar calculation shows that
the order of the stochastic wood $ \w_3 $
in Figure \ref{fig6} is 
$ \delta + \min(\delta,2\gamma) $
and that the order of 
the stochastic wood $ \w_4 $
in Figure \ref{fig7} 
is also $ \delta + \min(\delta,2\gamma) $.
Finally, we obtain that the stochastic
wood $ \w_5 $ in Figure \ref{fig8} is of order
$ \delta + 2 \min(\delta,\gamma) $.

\subsection{Trees and stochastic processes}

To identify each tree in $\textbf{ST}$ with a 
predictable stochastic process in $\mathcal{P}$, 
we define two functions
$\Phi$, $\Psi$ $:$ $\textbf{ST}$ $\rightarrow$ $\mathcal{P}$
recursively. For a given 
S-tree $\tr$ $=$ $(\tr',\tr'')$ 
with $k$ $=$ $\tr''( 1 )$ we 
define $\Phi( \tr )$ $:=$ $I^0_k$ 
when $ k = 0 $  
or $l( \tr )$ $=$ $1 $ and, 
when $ l( \tr )$ $\geq$ $2$ and $k$ $\neq$ $0$,
we define
$$
  \Phi( \tr )
  := I^n_k\left[
    \Phi( \tr_1 ), \dots, \Phi( \tr_n ) 
  \right],
$$
where $\tr_1$, $\dots$, $\tr_n$, $n$ $\geq$ $1$,  are the subtrees of $\tr$. 
In addition, for an arbitrary $ \tr $ in \textbf{ST}, we define $\Psi( \tr )$ $:=$ $0$ if $\tr$ is an active tree
and $\Psi( \tr )$ $=$ $\Phi( \tr )$ otherwise. Finally, for an 
S-wood $ \w$ $=$ $( \tr_1, \dots, \tr_n )$, 
$n \geq 1$, we define $ \Phi( \w )$ and $ \Psi( \w ) $
by
$$
  \Phi( \w ) = \Phi( \tr_1 ) + \ldots + \Phi( \tr_n ), \qquad
  \Psi( \w ) = \Psi( \tr_1 ) + \ldots + \Psi( \tr_n ) .
$$ 
As an example, for the elementary stochastic wood $ \w_0 $
(see Figure \ref{fig4}) we have
\begin{equation} \label{extr1a}
  \Phi( \w_0 ) = I^0_0 + I^0_{{1^*}} + I^0_{2^*} 
  \qquad
  \text{and}
  \qquad
  \Psi( \w_0 ) = I^0_0 .
\end{equation}
Hence, we obtain 
\begin{equation} \label{treesol1}
  \Phi( \w_0 ) = \Delta U 
\end{equation} 
from the equation above and equation \eqref{eqmain3}.
Since $(3,1)$ is an active node 
of $ \w_0 $, we obtain
\begin{equation} \label{extr1}
  \Phi( \w_1 ) 
  = I^0_0 + I^0_{1^*} + I^0_2 
  + I^1_{{2^*}}[ I^0_0 ]
  + I^1_{{2^*}}[ I^0_{{1^*}} ]
  + I^1_{{2^*}}[ I^0_{2^*} ]
\end{equation}
and 
\begin{equation} \label{extr2}
  \Psi( \w_1 ) 
  = I^0_0 + I^0_2 
\end{equation}
for the S-wood $ \w_1 = E_{(3,1)} \w_0 $  presented in Figure \ref{fig2}. 
Moreover, in view of equations \eqref{formel1}, \eqref{formel1b},
\eqref{formel2} and \eqref{formel2b}, we have
\begin{equation} \label{eformel}
  \Phi( \w ) 
  = \Phi( E_{(i,j)} \w )
\end{equation}
for every active node $ (i,j) $ of a 
stochastic wood $\w$ $\in$ $\textbf{SW}$.
Hence, we obtain 
$$
  \Phi( \w_1 ) = \Phi( E_{(3,1)} \w_0 )
  = \Phi( \w_0 )
$$
due to the equation above and the definition of $ \w_1 $.
Hence, we also obtain $ \Phi( \w_1 ) = \Delta U $.
And again equation \eqref{eformel} yields
$$
  \Phi( \w_5 ) = \Phi( \w_4 ) = \dots = \Phi( \w_0 )
  = \Delta U 
$$
by the definition of $ \w_5, \dots, \w_1 $.
As an further example, we obtain
\begin{equation} \label{eq01}
  \Phi( \w_2 ) = I^0_0 + I^0_1 + I^0_2 
  + I^1_{2^*}[ I^0_0 ]
  + I^1_{2^*}[ I^0_{1^*} ]
  + I^1_{2^*}[ I^0_{2^*} ]
  + I^1_{1^*}[ I^0_{0} ]
  + I^1_{1^*}[ I^0_{1^*} ]
  + I^1_{1^*}[ I^0_{2^*} ],
\end{equation}
which is, since 
$ \Phi( \w_2 ) = \Delta U $, 
nothing else than the expansion of the
solution given by equation \eqref{tay1rest}.
However, in the expression of the right hand 
side of the equation \eqref{eq01}, the solution
process $ ( U_t )_{ t \in [0,T] } $ occurs,
for example in the integrals 
$ I^1_{1^*}[ I^0_0 ] $.
More precisely, in every integral of type $ 1^* $
or $ 2^* $.
Therefore, we use $ \Psi( \w_2 ) $ instead
of $ \Phi( \w_2 ) $, which omits all these integrals,
as a good approximation of the solution process.
We have
$$
  \Delta U
  =
  \Phi( \w_2 ) 
  \approx
  \Psi( \w_2 ) = I^0_0 + I^0_1 + I^0_2 .
$$
We also note that the right hand side of 
the equation above is just the 
exponential Euler approximation
in equation \eqref{tay1restb}.
With the above notation and 
definitions we are now able 
to present the main result,
which is a presentation formula for the  
solution of the SPDE \eqref{eqspde} 
by the Taylor expansion
and an estimate of the remainder
of the Taylor approximation.
Hence, it can be seen that $ \Psi( \cdot ) $ is indeed
a good approximation of $ \Phi( \cdot ) $.
\begin{theorem} \label{mainthm}
Let Assumptions 
\ref{linearoperator}-\ref{initial}
be fulfilled.
Let ${\rm \mathbf{w} }$ be an arbitrary 
given stochastic wood in ${\rm \mathbf{SW}'}$
and let $ p \geq 1 $. 
Then there exists a constant 
$C > 0 $ such that 
\begin{equation} \label{taylorSPDE}
  \mathbb{P}\left[
  U_t
  = U_{ t_0 } + \Phi({\rm \mathbf{w} })(t)
  \right] = 1,
  \qquad
  \left| 
    U_t - \left( U_{ t_0 } + \Psi({\rm \mathbf{w} } )(t) \right)
  \right|_{ L^p }
  \leq C ( \Delta t )^{ \ord( {\rm \mathbf{w} } ) }
\end{equation}
for every $ t \in [t_0,T] $,
where $ U_t $, $ t \in [0,T] $, 
is the unique solution
of the SPDE \eqref{eqspde}.
Here 
the constant $ C > 0 $ 
only depends on 
the S-wood $ \w $,
on $ p \geq 1 $, on $ T > 0 $,
on $ \left| F(0) \right| $ given in 
Assumptions \ref{drift}, 
on $ \delta, ( L_i )_{ i \geq 0 } $
given in Assumptions \ref{stochconv}
and on $ ( R_i )_{ i \geq 1 }, 
( K_i )_{ i \geq 1 } $ given in
equation \eqref{constants}.
\end{theorem}
The representation of the solution here is 
a direct consequence of equations \eqref{treesol1} and 
\eqref{eformel} and the 
definition of $ \mathbf{SW}' $.
The proof for the estimate in \eqref{taylorSPDE}
will be  given in section \ref{secproofs}.
Here $ \Phi( \w ) = \Delta U $ is the 
Taylor expansion of the increment of the 
solution of the SPDE \eqref{eqspde}, while
$ \Psi( \w ) $ is the Taylor approximation 
of the increment of the solution 
and $\Phi( \w ) - \Psi( \w )$
is its  reminder,  which are  uniquely described by the S-wood $ \w $ in $ \textbf{SW}' $.
Since there are woods in $\textbf{SW}'$ with arbitrary high orders,  Taylor expansions  of  arbitrary high order 
can be constructed  by successively 
applying the operator $ E_{(i,j)}$  to 
the initial S-wood $ \w_0 $. 
Finally, the result of 
Theorem \ref{mainthm} can also 
be written as
$$
  U_t 
  = U_{ t_0 } 
  + \Psi( \w ) + O( ( \Delta t )^{ \ord(\w) } )
$$
for $ t \in [t_0,T] $ and a stochastic wood $ \w \in \textbf{SW}' $.

\section{Examples} \label{secex}

We present some examples here to illustrate 
the Taylor expansions introduced above 
and Assumptions 
\ref{linearoperator}-\ref{initial}.

\subsection{Abstract examples of the Taylor expansions} 
\label{secabsex}
We begin with some abstract examples 
of the Taylor expansions.

\subsubsection{Taylor expansion 
of order $\delta$} 
\label{sec511}
The first Taylor expansion of the solution is
given by the initial stochastic wood $ \w_0 $
(see Figure \ref{fig4}),
i.e. we have $ \Phi( \w_0 ) = \Delta U $ approximated
by $ \Psi( \w_0 ) $ with order $ \ord( \w_0 ) $.
Precisely, we have
$$
  \Psi( \w_0 )(t) 
  = ( e^{ A \Delta t } - I ) U_{ t_0 }
$$
and
$$
  \Phi( \w_0 )(t) 
  = ( e^{ A \Delta t } - I ) U_{ t_0 }
  + \int^t_{t_0} e^{ A(t-s) } F( U_s ) ds
  + \int^{t}_{t_0} e^{ A(t-s) } B( U_s ) dW_s
$$
almost surely due to equation \eqref{extr1a}.
Since, $ \ord( \w_0 ) = \delta $ 
(see section \ref{sectrees}.\ref{secorder}),
we finally obtain
$$
  U_t = e^{ A \Delta t } U_{ t_0 } 
  + O( ( \Delta t )^{ \delta } )
$$
as a Taylor expansion of order $ \delta $.
Here $ t $ is always in $ [t_0,T] $.

\subsubsection{Taylor expansion 
of order $ \delta + \min(\gamma,\delta) $} 
\label{sec512}
In the next step, we consider the Taylor expansion
given by the S-wood $ \w_1 $ (see Figure \ref{fig2}).
Here $ \Phi( \w_1 ) $ and $ \Psi( \w_1 ) $ are
presented in equation \eqref{extr1} and \eqref{extr2}.
Since $ \ord( \w_1 ) = \delta + \min( \gamma, \delta ) $
(see section \ref{sectrees}.\ref{secorder}), we obtain
$$
  U_t = 
    e^{ A \Delta t } U_{t_0}
    +
    \int^{t}_{t_0} e^{ A (t-r) } 
      B( U_{ t_0 } ) \,
    dW_r + O( (\Delta t)^{ 
      ( \delta + \min(\gamma,\delta) ) 
    } )
$$
as a Taylor expansion 
of order $ \delta + \min(\gamma,\delta) $.

\subsubsection{A further Taylor expansion 
of order $ \delta + \min(\gamma,\delta) $} 
\label{sec513}
The stochastic wood $ \w_2 $ (see Figure \ref{fig5}) has
order $ \delta + \min(\gamma,\delta) $ 
(see section \ref{sectrees}.\ref{secorder}).
Here the Taylor approximation 
$ \Psi( \w_2 ) $ of 
$ \Phi(\w_2) = \Delta U $
is given by
$
  \Psi( \w_2 )
  = I^0_0 + I^0_1 + I^0_2
$
and therefore, we obtain
\begin{eqnarray*}
  U_t &=& 
    e^{ A \Delta t } U_{t_0}
    + 
    A^{-1}(e^{ A \Delta t }-I) F( U_{t_0} )
    +
    \int^{t}_{t_0} e^{ A (t-r) } 
      B(U_{t_0})\,dW_r
    + O( (\Delta t)^{ 
      ( \delta + \min(\gamma, \delta) ) 
    } )
\end{eqnarray*}
as a further Taylor expansion 
of order $ \delta + \min( \gamma, \delta ) $.
This example corresponds in the case 
of additive
noise to the exponential Euler scheme
introduced in Jentzen \& Kloeden (2008b), 
which was already discussed 
in section \ref{sectay}.\ref{secsimple} 
(see also equation \eqref{onestep}).
One could ask, why one presents this Taylor expansion here,
if it is of the same order as the Taylor expansion above.
The reason is the following: The local
approximation quality of these both Taylor expansions
is the same. However, for numerical schemes
one has to consider the global approximation quality
of a one-step numerical scheme induced by a local
approximation. It turns out that the
one-step numerical scheme induced by this Taylor
expansion has very good global approximation properties,
which can also be seen in Jentzen \& Kloeden (2008b) 
in the case of
additive noise and is a task for further research in
the general case (see Remark \ref{fr}).

\subsubsection{Taylor expansion of order 
$ \delta + \min( 2 \gamma, \delta ) $} \label{sec514}
In the next step, the stochastic wood $ \w_3 $ 
(see Figure \ref{fig6}) has
order $ \delta + \min( 2 \gamma, \delta ) $ 
(see section \ref{sectrees}.\ref{secorder}).
Moreover, we have
$$
  \Psi( \w_3 )
  = I^0_0 + I^0_1 + I^0_2 
  + I^1_2[ I^0_0 ]
$$
and hence, we obtain
\begin{eqnarray*}
  U_t &=& 
    e^{ A \Delta t } U_{t_0}
    + 
    A^{-1}(e^{ A \Delta t }-I) F( U_{t_0} )
    +
    \int^{t}_{t_0} e^{ A (t-s) } 
      B(U_{t_0})\,dW_s
  \\
  &&
  + 
  \int^t_{t_0} e^{ A ( t-s) }
  B'( U_{t_0} )\left(
  \left( e^{ A \Delta s } - I \right) U_{ t_0} \right)
  dW_s 
    + O( (\Delta t)^{ 
      ( \delta + \min( 2 \gamma, \delta) ) 
    } )
\end{eqnarray*}
as a Taylor expansion of order 
$ \delta + \min( 2 \gamma, \delta ) $.

\subsubsection{Taylor expansion 
of order $ \delta + 2 \min( \gamma, \delta) $}
\label{sec515}
Since $ \w_4 $ yields the same Taylor approximation
as $ \w_3 $, i.e. $ \Psi( \w_4 ) = \Psi( \w_3 ) $,
we present here the Taylor approximation
for the S-wood $ \w_5 $ (see Figure \ref{fig8}).
We obtain
\begin{eqnarray*}
  U_t &=& 
    e^{ A \Delta t } U_{t_0}
    + 
    A^{-1}(e^{ A \Delta t }-I) F( U_{t_0} )
    +
    \int^{t}_{t_0} e^{ A (t-s) } 
      B(U_{t_0})\,dW_s
  \\
  &&
  + 
  \int^t_{t_0} e^{ A ( t-s) }
  B'( U_{t_0} )
  \int^s_{ t_0 } e^{ A ( s-r ) }
  B( U_{ t_0 } ) \,dW_r \, dW_s 
  \\
  &&
  + 
  \int^t_{t_0} e^{ A ( t-s) }
  B'( U_{t_0} )\left(
  \left( e^{ A \Delta s } - I \right) U_{ t_0} \right)
  dW_s 
    + O( (\Delta t)^{ 
      ( \delta + 2\min(\gamma, \delta) ) 
    } )
\end{eqnarray*}
as a Taylor expansion of order 
$ \delta + 2 \min( \gamma, \delta) $
(see section \ref{sectrees}.\ref{secorder}).
This example corresponds to the Taylor expansion
introduced in section \ref{sectay}.\ref{sechigher}.

\subsection{Additive noise}

In the case of additive noise the stochastic trees
and stochastic woods can be simplified.
Therefore, a detailed presentation of abstract 
and concrete
examples of the Taylor expansions in the case
of additive noise can be found in 
Jentzen \& Kloeden (2008c).

\subsection{Noise via multiplication operators}

As an example of our assumptions, we 
consider the stochastic heat equation
with multiplicative space time noise in
one dimension
(see also section 2 in
M\"uller-Gronbach \& Ritter (2007a) 
for a very precise description 
of this example).
Let $ H = L^2\left( (0,1), \mathbb{R} \right) $ be the
Hilbert space of all square integrable functions from 
$ (0,1) $ to $ \mathbb{R} $.
The scalar product and the norm in $ H $ are given by
$$
  \left< f, g \right> = \int^1_{ 0 } f(x) g(x) \, dx, \qquad
  \left| f \right| = \left( \int^1_{ 0 } f(x)^2 \, dx 
  \right)^{ \frac{1}{2} }
$$
for every $ f, g \in H $. 
We also define $ U := H $.
Moreover, let $ A = \Delta $
be the Laplacian with Dirichlet boundary conditions
in one dimension, i.e.
$ \mathcal{I} = \mathbb{N} $ and 
$$ 
  e_i(x) = \sqrt{2}
  \sin( i \pi x ), \qquad
  \lambda_i =
  \pi^2 i^2
$$
for all
$ x \in (0,1) $
and all $ i \in \mathbb{N} $.
The operator $ A $ is then given by
$$
  A f = \sum_{ i \in \mathcal{I} } 
  - \lambda_i
  \left< e_i, f \right> e_i
$$
for all $ f \in D(A) $
with 
$ D(A) = 
\left\{ f \in H \; \big| \;
\sum_{ i \in \mathcal{I} } 
\lambda_i^2
\left| \left< e_i, f \right> \right|^2 < \infty
\right\} $
and in that way Assumption \ref{linearoperator} is fulfilled.
Furthermore let $ F = 0 $ and we define $ B $ by
\begin{equation} \label{defB}
  B : H \rightarrow L(H,D),
  \qquad
  \left( B(v)(w) \right)(x) 
  := v(x) \cdot w(x)
\end{equation}
for every
$ x \in (0,1) $ and $ v, w \in H $,
where $ D = L^1\left( (0,1), \mathbb{R} \right) $
is the Banach space of all integrable functions from
$ (0,1) $ to $ \mathbb{R} $.
Indeed, $ B $ is well defined, since 
\begin{eqnarray*}
  \left| B(v)(w) \right|_{D}
  & = &
  \int^1_0 \left| v(x) \cdot w(x) \right| dx
  \leq
  \left( 
    \int^1_0 \left| v(x) \right|^2 dx 
  \right)^{ \frac{1}{2} }
  \left( 
    \int^1_0 \left| w(x) \right|^2 dx 
  \right)^{ \frac{1}{2} }
  \\
  & = &
  \left| v \right| \cdot \left| w \right|
\end{eqnarray*}
for all $ v, w \in H $ 
due to the inequality of Cauchy-Schwartz 
and therefore
$ B(v) $ is indeed a bounded linear operator 
from $ H $ to $ D $
with 
$
  \left| B(v) \right|_{ L(H,D) } 
  \leq \left| v \right|
$
for all $ v \in H $.
Even more, by definition, also $ B $
is a linear operator from $ H $ to
$ L(H,D) $.
In particular, it is infinitely
often differentiable with the
derivatives $ B'(v) = B $ and 
$ B^{(i)}(v) = 0 $ for all $ i \geq 2 $ and 
every $ v \in H $.
Furthermore, it is straightforward to verify that
$ (-A)^\gamma e^{ A t } B(v) $ 
is a Hilbert-Schmidt operator 
from $ H $ to $ H $ with 
$$
  \left|
  e^{ A t } B( v )
  \right|_{HS}
  \leq
  4 (T+1)
  \left| v \right|
  t^{ - \frac{1}{4} } ,
  \qquad
  \left|
  (-A)^\gamma e^{ A t } B(v)
  \right|_{ HS }
  \leq 4 (T+1)
  \left| v \right|
  t^{ -(\frac{1}{4} + \gamma) }
$$
for every $ v \in H $, $ t \in (0,T] $, 
$ \gamma \in (0, \frac{1}{4}) $
(see Jentzen \& Kloeden (2008d) 
and also Remark 2 in 
M\"uller-Gronbach \& Ritter (2007a)).
Hence,
since $ L^1(0,1) \subset D( (-A)^{(-\frac{1}{2})} )$
continuously, 
Assumptions \ref{drift} and \ref{stochconv} are also fulfilled
with the parameters
$
  \gamma = \frac{1}{4} - r
$
and
$
  \delta = \frac{1}{4}
$
for an arbitrary but very small $ r > 0 $.
Now, we present the abstract Taylor expansions
in section \ref{secex}.\ref{secabsex} in that
particular situation.
First of all, the Taylor expansion 
in section \ref{secex}.\ref{secabsex}.\ref{sec511}
is here given by
$$
  U_t = e^{ A \Delta t } U_{ t_0 } 
  + O( ( \Delta t )^{ \frac{1}{4} } ) .
$$
Here $ t $ is always in $ [t_0,T] $.
In the next step, we consider the Taylor expansion
in section \ref{secex}.\ref{secabsex}.\ref{sec512}
respectively 
\ref{secex}.\ref{secabsex}.\ref{sec513},
which is here given by
$$
  U_t = 
    e^{ A \Delta t } U_{t_0}
    +
    \int^{t}_{t_0} e^{ A (t-r) } 
      B( U_{ t_0 } ) \,
    dW_r + O( (\Delta t)^{ 
      ( \frac{1}{2} - r ) 
    } ) .
$$
Moreover, we have
\begin{eqnarray*}
  U_t 
  &=& 
    e^{ A \Delta t } U_{t_0}
    +
    \int^{t}_{t_0} e^{ A (t-s) } 
      B(U_{t_0})\,dW_s
  \\
  &&
  + 
  \int^t_{t_0} e^{ A ( t-s) }
  B\left(
  \left( e^{ A \Delta s } - I \right) U_{ t_0} \right)
  dW_s 
    + O( (\Delta t)^{ 
      \frac{1}{2} 
    } )
  \\
  &=& 
    e^{ A \Delta t } U_{t_0}
    +
    \int^{t}_{t_0} e^{ A (t-s) } 
      B( e^{ A \Delta s } U_{t_0})\,dW_s
    + O( (\Delta t)^{ 
      \frac{1}{2} 
    } )
\end{eqnarray*}
for the Taylor expansion in section
\ref{secex}.\ref{secabsex}.\ref{sec514}.
Finally, we obtain
\begin{eqnarray*}
  U_t &=& 
    e^{ A \Delta t } U_{t_0}
    + 
    \int^{t}_{t_0} e^{ A (t-s) } 
      B( e^{ A \Delta s } U_{t_0}) \, dW_s
  \\
  &&
  + 
  \int^t_{t_0} e^{ A ( t-s) }
  B\left( \int^s_{ t_0 } e^{ A ( s-r ) }
  B( U_{ t_0 } ) \,dW_r \right) dW_s 
    + O( (\Delta t)^{ 
      ( \frac{3}{4} - 2 r ) 
    } )
\end{eqnarray*}
for the Taylor expansion in section
\ref{secex}.\ref{secabsex}.\ref{sec515}.
\begin{remark} \label{fr}
For further research, 
it would be interesting to 
analyze the numerical one-step scheme induced
by the Taylor expansion in section 
\ref{secex}.\ref{secabsex}.\ref{sec513}.
Let $ 0 = t_0 < t_1 < \dots < t_N = T $
be given discretization times.
Then, one could consider
the scheme $ ( X_n )_{ n \geq 1 } $ 
given by
$
  X_0 = u_0
$
and
$$
  X_{n+1} = 
    e^{ A (t_{n+1} - t_n) } X_n
    + 
    \int^{t_{n+1}}_{t_n} e^{ A (t_{n+1}-r) } 
      B(X_n)\,dW_r
$$
for $ n = 0, 1, \dots, N-1 $.
Of course, this is just a discretization in time,
but for the space discretization one could use
standard methods 
such as spectral Galerkin methods.
Note also that the second summand in the scheme
given above is an It$\hat{o}$ integral with an
integrand that is $ \mathcal{F}_{ t_n } $-measurable
and therefore the conditional distribution 
of this integral
(with respect to $ \mathcal{F}_{ t_n } $)
is normal.
Although it was shown 
in the seminal paper
Davie \& Gaines (2000) that 
even for a large class of 
linear functionals
the computational order barrier 
of $ \frac{1}{6} $ (convergence speed
of the linear implicit Euler scheme)
cannot be overcome
in this special example,
there is still hope that this
numerical scheme converges with
a higher order, since the linear functionals
above do not fit in the class
of linear functionals used there.
Therefore, from my point of view, it 
would be very interesting
to analyze the convergence order of the
scheme given above.
\end{remark}

\section{Proofs}\label{secproofs}

In the following,
we need a version of 
the Burkholder-Davis-Gundy inequality
in infinite dimensions (see Lemma 7.7 in 
Da Prato \& Zabczyk (1992)).
\begin{lemma} \label{lem2}
Let $ X_t $, $ t \in [t_0,T] $ be
a predictable stochastic process, whose
values are Hilbert-Schmidt operators
from $ U $ to $ H $ with 
$ \mathbb{E} 
\int^T_{t_0} \left| X_s \right|^2_{HS} ds < \infty $. 
Then, 
$$
  \left| 
      \int^t_{t_0} X_s \, dW_s
  \right|_{ L^p }
  \leq
  p
  \left(
    \int^t_{t_0}
      \left|
      \left| 
        X_s 
      \right|_{ HS }  
      \right|_{ L^p }^2
    ds
  \right)^{ \frac{1}{2} }
$$
for every $ t \in [t_0,T] $ and every $ p \geq 2 $.
(Both sides could be infinite.)
\end{lemma}

\subsection{Proof of Theorem \ref{mainthm}}

In view of the  definitions of the mappings $ \Phi $ and $ \Psi $,  
the proof of the Theorem \ref{mainthm} immediately follows from the next lemma.
\begin{lemma}
Let Assumptions 
\ref{linearoperator}-\ref{initial}
be fulfilled.
Let $\tr$ $=$ $( \tr', \tr'')$ be a given S-tree
and let $ p \geq 1 $. 
Then there exists a constant 
$C > 0 $ such that 
$$
  \left| \Phi( \tr )(s) \right|_{ L^p }
  \leq C \cdot ( \Delta s )^{ \ord( \tr ) }
$$
for all $ s \in [t_0,T] $, 
where $C>0$ only depends on 
the S-tree $ \tr $,
on $ p \geq 1 $,
on $ \left| F(0) \right| $ given in 
Assumptions \ref{drift}, 
on $ \delta, ( L_i )_{ i \geq 1 } $
given in Assumptions \ref{stochconv}
and on $ ( R )_{ i \geq 1 }, 
( K_i )_{ i \geq 1 } $ given in
equation \eqref{constants}.
\end{lemma}
\begin{proof} 
We can assume without loss of generality 
that $ p\geq 2 $.
Throughout the proof $ C > 0 $ is a constant
changing from line to line but only depends
on $ \left| F(0) \right| $, $ \delta $, 
$( R_i )_{ i \geq 1 }$, $p$, $\tr$, 
$( L_i )_{ i \geq 1 } $ and 
$ ( K_i )_{ i \geq 1 } $.
We will prove 
the  assertion by induction 
with respect to the number of nodes $l(\tr )$. \\
Let $ s \in [t_0,T] $ and 
suppose that $k = \tr''(1)$. 
Then, by Lemma \ref{lem2}, we have
\begin{eqnarray*}
  \left|
    \left( 
      e^{ A \Delta s } - I 
    \right) U_{ t_0 }
  \right|_{ L^p }
  & \leq &
  \left|
    (-A)^{-\gamma}
    \left( 
      e^{ A \Delta s } - I 
    \right) 
  \right|
  \cdot  
  \left|
    (-A)^{\gamma} U_{ t_0 }
  \right|_{ L^p } \\[2ex]
  & \leq &
  R_p \, ( \Delta s )^{ \gamma } \leq C ( \Delta s )^{\gamma} ,
  \\
  \left| 
    \int^{s}_{t_0} e^{ A(s-r) } F( U_r ) \, dr
  \right|_{ L^p }
  & \leq &
  \int^{s}_{t_0}
    \left|
      e^{ A(s-r) } F( U_r )
    \right|_{ L^p } \, 
  dr
  \leq
    \int^{s}_{t_0}
      \left( \left| F(0) \right| +
        K_0
        \left|
          U_r
        \right|_{ L^p }  
      \right)\,
    dr
  \\ 
  & \leq &
      \left( \left| F(0) \right| +
        K_0 R_p
      \right) ( \Delta s )
  \leq C ( \Delta s )
\end{eqnarray*}
and
\begin{eqnarray*}
  \left| 
  \int^{s}_{t_0}  e^{ A(s-r) } B( U_r ) \, dW_r 
  \right|_{ L^p }
  & \leq &
  p
  \left(
  \int^{s}_{t_0} \left| 
  \left|e^{ A(s-r) } B( U_r ) \right|_{HS} 
  \right|_{ L^p}^2 \,
    dr \right)^{ \frac{1}{2} } \\[2ex]
  & \leq &
    L_0 \, p
    \left(
    \int^{s}_{t_0}
        \left( 1 + \left| U_r \right|_{L^p} \right)^2 
        ( s - r )^{ (2 \delta - 1) } 
    dr \right)^{ \frac{1}{2} }
  \\
  & \leq &
    ( 1 + R_p ) L_0 \, p 
    \left( 
    \int^{ s }_{ t_0 } (s-r)^{ (2 \delta - 1) } dr 
    \right)^{ \frac{1}{2} }
  \\
  &=&
    ( 1 + R_p ) L_0 \, p 
    \left( 
      \int^{ \Delta s }_{ 0 } r^{ (2 \delta - 1) } dr
    \right)^{ \frac{1}{2} }
  \\
  &=& 
    \left(
    \frac{ ( 1 + R_p ) L_0 \, p }{ \sqrt{ 2 \delta } }
    \right) 
    ( \Delta s )^{ \delta } 
    \leq 
    C ( \Delta s )^{ \delta } .
\end{eqnarray*}
This yields
$
  \left| \Phi( \tr )(s) \right|_{ L^p } 
  \leq C ( \Delta s )^{ \ord( \tr ) } 
$
when $k=0$ or $l( \tr )$ $=$ $1 $
since by  definition $\Phi( \tr )$ $=$ $I^0_k $ here.\\
Suppose now that $l( \tr )$ $\geq$ $2$ and $ k \neq 0 $.
Let  $ \tr_1$, $\dots$, $\tr_n $, $ n \geq 1 $, be   the subtrees of $ \tr $. 
Then, by definition, we have
$$
\Phi( \tr )(s) = I^n_k\left[   \Phi( \tr_1 ), \dots, \Phi( \tr_n )  \right](s)
$$
for $ s \in [t_0,T] $. 
When $ k = {1^*} $, we have
$$
  \Phi( \tr )(s) 
  =
  \int^{s}_{t_0}  e^{ A(s-r) } 
  \left(  
    \int^1_0 F^{(n)} 
    ( U_{ t_0 } + \theta \Delta U_r )
    \left( \Phi( \tr_1 )(r), \dots, \Phi( \tr_n )(r) \right) 
    \frac{ (1-\theta)^{(n-1)} }{ (n-1)! } \, d\theta \right)
  \,dr
$$
and therefore
$$
  \left| \Phi( \tr )(s) \right|
  \leq 
  K_n
  \int^{s}_{t_0}
    \left| \Phi( \tr_1 )(r) \right| \cdots 
    \left| \Phi( \tr_n )(r) \right| \, dr .
$$
Hence, we obtain
\begin{eqnarray*}
\left| \Phi( \tr )(s) \right|_{ L^p }  
& \leq &
  K_n
  \int^{s}_{t_0}\left| \left| \Phi( \tr_1 )(r) \right| 
  \cdots 
  \left| \Phi( \tr_n )(r) \right|\, 
  \right|_{ L^p } \,  dr 
  \\[2ex]
  & \leq &
  K_n
  \int^{s}_{t_0} 
  \left|  \Phi( \tr_1 )(r)  \right|_{ L^{p n } } 
  \cdots
  \left| \Phi( \tr_n )(r) \right|_{ L^{p n } } \, dr 
  \\[2ex]
  & \leq &
  C  ( \Delta s )^{ ( 1 + \ord(\tr_1) + \ldots + \ord(\tr_n) ) }
  = C  ( \Delta s )^{ \ord( \tr ) },
\end{eqnarray*}
since $ l(\tr_1)$, $\dots$, $l(\tr_n)$ $\leq$ $l(\tr)-1 $ and
we can apply the induction assertion to the subtrees. 
A similar calculation shows 
the result when $ k = 1 $. \\
When $ k = {2^*} $, we have
$$
\Phi( \tr )(s)  =   \int^{s}_{t_0}  
  e^{ A(s-r) }  
  \left( 
    \int^1_0 B^{(n)}( U_{ t_0 } + \theta \Delta U_r )
    \left( \Phi( \tr_1 )(r), \dots, \Phi( \tr_n )(r) \right) 
    \frac{ (1-\theta)^{(n-1)} }{ (n-1)! } 
  d\theta \right) dW_r
$$
and therefore
\begin{eqnarray*}
  &&
  \left| \Phi( \tr )(s) \right|_{ L^p }
  \\
  & \leq &
  p 
  \left(
  \int^{s}_{t_0}
  \left|
  \left|
      \int^1_0 
      e^{ A(s-r) }  
      B^{(n)}( U_{ t_0 } + \theta \Delta U_r )
      \left( \Phi( \tr_1 )(r), \dots, \Phi( \tr_n )(r) \right)
      \frac{ (1-\theta)^{(n-1)} }{ (n-1)! } 
    d\theta
  \right|_{HS} 
  \right|_{ L^p }^2 \, dr
  \right)^{ \frac{1}{2} } 
  \\
  & \leq &
  p 
  \left(
  \int^{s}_{t_0}
    \left(
    \int^1_0 
      \left|
      \left|
      e^{ A(s-r) }  
      B^{(n)}( U_{ t_0 } + \theta \Delta U_r )
      \left( \Phi( \tr_1 )(r), \dots, \Phi( \tr_n )(r) \right) \, 
      \right|_{HS} 
      \right|_{ L^p }
    d\theta 
  \right)^2 \, dr
  \right)^{ \frac{1}{2} } 
  \\
  & \leq &
  L_n p 
  \left(
  \int^{s}_{t_0}
    (s - r)^{ (2 \delta - 1) }
    \left(
    \int^1_0 
      \left|
      \left( 1 + \left| U_{ t_0 } + \theta \Delta U_r \right| \right)
      \left| \Phi( \tr_1 )(r) \right| 
      \dots \left| \Phi( \tr_n )(r) \right| 
      \right|_{ L^p }
    d\theta 
  \right)^2 \, dr
  \right)^{ \frac{1}{2} } .
\end{eqnarray*}
Since
\begin{eqnarray*}
  &&
    \int^1_0 
      \left|
      \left( 
        1  
        + \left| U_{ t_0 } + \theta \Delta U_r \right| 
      \right)
      \left| \Phi( \tr_1 )(r) \right| 
      \dots \left| \Phi( \tr_n )(r) \right| 
      \right|_{ L^p }
    d\theta 
  \\ 
  & \leq &
    \int^1_0 
      \left( 
        1  
        + \left| U_{ t_0 } + \theta \Delta U_r \right|_{ L^{p(n+1)} }
      \right)
        \left| 
          \Phi( \tr_1 )(r) 
        \right|_{ L^{p(n+1)} } 
        \cdots 
        \left| 
          \Phi( \tr_n )(r)
        \right|_{ L^{ p(n+1) } }
    d\theta 
  \\
  & \leq &
    C \, 
    (r-t_0)^{ 
      \left( 
        \ord( \tr_1 ) + \ldots + \ord( \tr_n )
      \right)
    }
    \int^1_0 
        \left( 1 +
        \left| 
          \theta U_r + (1-\theta) U_{ t_0 }
        \right|_{ L^{ p(n+1) } }
        \right)
    d\theta 
  \\
  & \leq &
    C \, \left( 1 + R_{ p(n+1) } \right) \,
    ( r - t_0 )^{ 
      \left( 
        \ord( \tr_1 ) + \ldots + \ord( \tr_n )
      \right)
    }
  \leq
    C
    ( \Delta s )^{ 
      \left( 
        \ord( \tr_1 ) + \ldots + \ord( \tr_n )
      \right)
    } 
\end{eqnarray*}
for every $ r \in [t_0,s] $
due to the induction assertion, we obtain
\begin{eqnarray*}
  \left| \Phi( \tr )(s) \right|_{ L^p }  
  & \leq &
  C
  \left(
  \int^{s}_{t_0}
    (s - r)^{ (2 \delta - 1) }
  \, dr
  \right)^{ \frac{1}{2} }
  ( \Delta s )^{ \left( 
    \ord( \tr_1 ) + \dots + \ord( \tr_n )
  \right) }
  \\
  & \leq &
  C
  ( \Delta s )^{ \left( \delta + 
    \ord( \tr_1 ) + \dots + \ord( \tr_n )
  \right) }
  =
  C
  ( \Delta s )^{ \ord( \tr ) } .
\end{eqnarray*}
A similar calculation shows the result when $ k = 2 $. 
\end{proof}

\begin{acknowledgements}
This work was supported by 
the DFG project ``Pathwise numerics and dynamics
of stochastic evolution equations''.
\end{acknowledgements}

\end{document}